\documentclass[%
  onecolumn
   , colorlinks 
]{mpi2015-cscpreprint}

\usepackage[american]{babel}

\usepackage{graphicx}

\usepackage{amssymb}
\usepackage{amsthm}

\newcommand{\pde}{\texttt{PDE}}
\newcommand{\pinn}{\texttt{PINN}}
\newcommand{\fdm}{\texttt{FDM}}
\newcommand{\fem}{\texttt{FEM}}
\newcommand{\dnn}{\texttt{DNN}}
\newcommand{\fde}{\texttt{FDE}}
\newcommand{\sde}{\texttt{SDE}}
\newcommand{\sindy}{\texttt{SINDy}}
\newcommand{\propod}{\texttt{POD}}
\newcommand{\nn}{\texttt{NN}}
\newcommand{\deim}{\texttt{DEIM}}
\newcommand{\eim}{\texttt{EIM}}
\newcommand{\qdeim}{\texttt{Q-DEIM}}
\newcommand{\maxm}{\texttt{max}}
\newcommand{\minm}{\texttt{min}}
\newcommand{\SVD}{\texttt{SVD}}

\newcommand{\deepnn}{\texttt{DNN}}
\newcommand{\gspinn}{\texttt{GS-PINN}}
\newcommand{\rspinn}{\texttt{RS-PINN}}

\usepackage[vlined,longend,linesnumbered,ruled]{algorithm2e}

\usepackage[american]{babel}

\usepackage{graphicx}

\usepackage{amssymb}
\usepackage{amsthm}
\usepackage{amsmath}

\usepackage{caption}
\usepackage{subcaption}

\usepackage{graphics} 
\usepackage{epsfig} 
\usepackage{tikz}
\usepackage{xcolor}
\usepackage[latin1]{inputenc}
\usepackage[T1]{fontenc}

\usepackage{cleveref}

\usepackage{multirow}

\usepackage{comment}

\usepackage{longtable,array,amsmath,booktabs,tabularx,ragged2e}

\usetikzlibrary{arrows,automata, quotes, calc, trees, positioning, chains, shapes.geometric, decorations.pathreplacing, decorations.pathmorphing, shapes, matrix, shapes.symbols}

\usepackage{circuitikz}

\usetikzlibrary{arrows}
\usepackage[noend]{algpseudocode}

\usepackage{lscape}
\usepackage[colorinlistoftodos,prependcaption,textsize=footnotesize]{todonotes}

\usepackage{hyperref}
\usepackage{enumitem}

\usepackage[software, hardware]{mymacros}

\begin{document}
  

\title{GS-PINN: Greedy Sampling for Parameter Estimation in Partial Differential Equations}
  
\author[$\ast$]{Ali Forootani}
\affil[$\ast$]{Max Planck Institute for Dynamics of Complex Technical Systems, 39106 Magdeburg, Germany.\authorcr%
  \email{forootani@mpi-magdeburg.mpg.de}, \orcid{0000-0001-7612-4016}}
  
\author[$\ast$]{Harshit Kapadia}
\affil{ \email{kapadia@mpi-magdeburg.mpg.de }, \orcid{0000-0003-3214-0713}}

\author[$\ast$]{Sridhar Chellappa}
\affil{
  \email{chellappa@mpi-magdeburg.mpg.de}, \orcid{0000-0002-7288-3880}}

\author[$\ast$]{Pawan Goyal}
\affil{ \email{goyalp@mpi-magdeburg.mpg.de}, \orcid{0000-0003-3072-7780} }

\author[$\ast$]{Peter Benner}
\affil{ \email{benner@mpi-magdeburg.mpg.de}, \orcid{0000-0003-3362-4103} }

\shorttitle{Example short title}
\shortauthor{F. Author, S. Author}
\shortdate{}
  
\keywords{Discrete Empirical Interpolation Method (\texttt{DEIM}), Physics Informed Neural Network (\pinn), Parameter Estimation}

\msc{MSC1, MSC2, MSC3}
  
\abstract{
Partial differential equation parameter estimation is a mathematical and computational process used to estimate the unknown parameters in a partial differential equation model from observational data. This paper employs a greedy sampling approach based on the Discrete Empirical Interpolation Method to identify the most informative samples in a dataset associated with a partial differential equation to estimate its parameters. Greedy samples are used to train a physics-informed neural network architecture which maps the nonlinear relation between spatio-temporal data and the measured values. To prove the impact of greedy samples on the training of the physics-informed neural network for parameter estimation of a partial differential equation, their performance is compared with random samples taken from the given dataset. Our simulation results show that for all considered partial differential equations, greedy samples outperform random samples, i.e., we can estimate parameters with a significantly lower number of samples while simultaneously reducing the relative estimation error. A \texttt{Python} package is also prepared to support different phases of the proposed algorithm, including data prepossessing, greedy sampling, neural network training, and comparison.
  
}

\novelty {Within this study, we employ the discrete empirical interpolation method (\deim), specifically the QR-factorization-based variant known as \qdeim, as a strategic sampling technique for mitigating the computational complexities and time demands associated with parameter estimation for partial differential equations (\pde s) using neural networks. Our methodology involves the judicious pre-selection of spatio-temporal data, thereby constructing a reduced dataset for training a neural network to estimate the coefficients of the underlying \pde~governing the data. We establish that our proposed \qdeim-based sampling approach not only reduces the required training data for the neural network but also yields a commendable approximation of \pde~coefficients in fewer training iterations.}

\maketitle

  
\section{Introduction}%
\label{sec:intro}
Partial differential equations (\pde s) are fundamental mathematical tools used to describe the behavior of physical systems in various scientific and engineering disciplines, such as fluid dynamics, heat transfer, and quantum mechanics. Accurately identifying the underlying \pde\ governing a given system is a crucial step in understanding its behavior and making predictions. Traditional methods for identifying \pde s often rely on domain knowledge, mathematical derivations, and experimental data, which can be labor-intensive and may not be applicable in complex or poorly understood systems.

In the literature, most attempts have been made to solve the \pde s analytically or numerically. Analytical methods are based on finding a change of variable to transform the equation into something soluble or on finding an integral form of the solution \cite{de2006integral}. These methods are often used for simple \pde s, but they can be difficult to apply to more complex equations. Numerical methods approximate the solution of a \pde~by discretizing the domain of the equation and solving a system of algebraic equations, such as the finite difference method (\fdm)~\cite{zeneli2021numerical} and the finite element method (\fem)~\cite{bai2022local}. These methods are more general than their analytical counterparts, but they can be computationally expensive.

Besides solving the \pde s with conventional methods, recent advances in machine learning techniques have proved their potential to address \pde~problems in scenarios with limited data. This implies having access solely to the \pde~problem data, rather than an extensive set of value pairs for the independent and dependent variables \cite{blechschmidt2021three}. Taking advantage of modern machine learning software environments has provided automatic differentiation capabilities for functions realized by deep neural networks (\dnn) which is a mesh-free approach and can break the curse of dimensionality of the conventional methods~\cite{lu2021deepxde}. This approach was introduced in \cite{raissi2019physics}, where the term physics-informed neural networks~(\pinn s) was coined.

\pinn s combine the expressive power of neural networks with the physical knowledge of governing equations. \pinn s can be used to solve a wide variety of problems, including (i) Forward problems: to solve nonlinear \pde s on arbitrary complex-geometry domains \cite{chen2019physics, raissi2020hidden}, including fractional differential equations (\fde s)~\cite{pang2019fpinns}, and stochastic differential equations~(\sde s~\cite{zhang2019quantifying,yang2018physics,nabian2018deep,zhang2019learning}, (ii) Inverse problems: to retrieve the unknown parameters in the \pde~\cite{raissi2019physics, chen2019physics}, (iii) Control problems: designing control inputs to achieve a desired state~\cite{antonelo2021physics}, and (iv) Optimization problems: finding the optimal solution to a problem subject to constraints~\cite{mowlavi2023optimal}.

Unlike the previous literature that focused more on solving \pde s, in this work we consider solving the inverse problem, that is, obtaining the underlying \pde s based solely on time series data. One of the advantages of using \pinn~for inverse problem is that it requires minimum modifications of the Deep Neural Network (\dnn) architecture to recover the \pde~ coefficients and underlying governing equations of the given time series dataset.

When it comes to the challenge of discovering the underlying \pde~governing a given set of spatio-temporal data, choosing an appropriate collocation/sampling strategy is crucial. As the amount of data increases, there is a need for sampling strategies that require fewer collocation points to recover the underlying dynamical system. Indeed, the location and distribution of these sampling points are highly important to the performance of \pinn s. In literature, a few simple residual point sampling methods have mainly been employed, among them we can report non-adaptive uniform sampling such as (i) equispaced uniform grid, (ii) uniform random sampling, (iii) Latin hypercube sampling \cite{raissi2019physics}, (iv) Halton sequence \cite{halton1960efficiency}, (v) Hammersley sequence~\cite{Wongetal97}, (vi) Sobol sequence \cite{sobol1967distribution}; and adaptive non-uniform sampling such as residual-based adaptive distribution (\texttt{RAD}) and residual-based adaptive refinement with distribution (\texttt{RAR-D}) \cite{wu2023comprehensive}. Even though these methods seem promising, they are highly problem-dependent and are usually tedious and time-consuming. In particular, the case of adaptive sampling such as \texttt{RAD} or \texttt{RAR-D} requires re-sampling of the dataset within the training loop of the Neural Network (\nn), which adds significant computational cost to the corresponding algorithm. For further details regarding non-adaptive and residual-based adaptive sampling strategies and their comparisons in \pinn~training, we refer to the work reported in \cite{wu2023comprehensive}.

In this paper, we make use of the discrete empirical interpolation method (\deim) as a sampling technique to reduce the computational complexity and time consumption of \pde~parameter estimation using neural networks. In a nutshell, our work consists of (a) pre-selecting a portion of the available spatio-temporal data for a \pde~and (b) training a neural network using this reduced dataset in order to estimate the coefficients of the underlying \pde~governing the data. We demonstrate that our proposed sampling approach not only reduces the data needed to train the neural network but also recovers a good approximation of the \pde~coefficients in fewer iterations of training.

The approach we use to perform the sampling, viz., \deim~was originally proposed in the context of model order reduction to significantly reduce the computational complexity of the popular \propod~method for constructing reduced-order models for time-dependent and/or parametrized nonlinear \pde s \cite{chaturantabut2010nonlinear}. It is worth highlighting that \deim~is a discrete variant of the empirical interpolation method (\eim) for constructing an approximation of a non-affine parameterized function with a spatial variable defined in a continuous bounded domain with associated error bound on the quality of approximation \cite{barrault2004empirical}. In this paper, we particularly employ the QR-factorization-based \deim~procedure, i.e., \qdeim, which has a better upper bound error with respect to the original \deim~algorithm and has numerically robust high-performance procedures, already available in software packages such as \texttt{Python}, \texttt{LAPACK}, \texttt{ScaLAPACK}, and \texttt{MATLAB} \cite{drmac2016new}. More specifically, we methodically investigate the impact of the most informative samples acquired via \qdeim~on \pde~snapshot matrix for estimating the parameters associated with the \pde. We also compare the result of \qdeim~sampling in \pde~parameter estimation with the ones that we computed with the random sampling approach. Our findings prove the significant impact of most informative samples on the training of \pinn~architecture for the \pde~parameter estimation. In our setting, the \pinn~architecture which takes advantage of \qdeim~is named greedy sampling based \pinn~(\gspinn) while in contrast the one that is fed with random samples will be called random sampling based \pinn~(\rspinn). A \texttt{Python} package is prepared to support different implementation phases of \gspinn~and its comparison with \rspinn~on estimating the \pde~parameters corresponding to the Allen-Cahn equation, Burgers' equation and  Korteweg-de Vries equation.

This article is organised as follows:
In \Cref{sec:background}, an overview of \pde~parameter estimation will be discussed. We explain our greedy sampling approach and its usage for selecting most informative samples from the dataset corresponding to \pde s in \Cref{sec:qdeim}.
In \Cref{sec:method} we introduce the core algorithm for the \pde~parameter estimation. \Cref{sec:Simulation} is devoted to the simulation and comparison of the \gspinn~ with \rspinn. The paper is concluded in \Cref{sec:conclusion}.

\section{An Overview of \pde~parameter estimation}%
\label{sec:background}
In this section, we provide a brief background regarding the \pde~parameter estimation problem. Our aim is not to provide an in-depth discussion of this concept but rather to present a brief overview of the concepts narratively.

\subsection{Recap on \pde~parameters estimation}


The estimation of parameters for \pde s holds significance in various domains, including geophysical exploration and medical images \cite{haber2014computational}. This process is often framed as an optimization challenge constrained by \pde s, typically addressed iteratively through gradient-based optimization techniques \cite{fung2019multiscale}. Although significant progress has been made over the last two decades involving high-order schemes for \pde s, automatic differentiation (AD), \pde-constrained optimization, and optimization under uncertainty, parameter estimation in large-scale problem remains a significant challenge \cite{tartakovsky2020physics}.

Consider a nonlinear \pde~of the form
\begin{equation}\label{pde_1}
	\frac{\partial \bu}{\partial t} = \mathbf{f}_p(\mathbf{u},\mathbf{u}_x, \mathbf{u}_{xx}, \dots, x),
\end{equation}
where $t \in [0,t_{\maxm} ]$ is the time variable, $x \in [x_{\minm}, x_{\maxm}] $ is the space variable, $\mathbf{u}(x, t)$ is the solution, and $\mathbf{f}_p$ is generally a non-linear function of the solution and its derivatives. Since the function $\mathbf{f}_p(\cdot)$ depends on the parameter vector $\mathbf{p} \in \R^k$, the solution $\bu(\cdot)$ is also a function of $\mathbf{p}$. Moreover, in this work, we consider that the parameter vector $\mathbf{p}$ is a priori unknown, and we seek to estimate it.

The main assumption behind our parameter estimation is that the function $\mathbf{f}_p(\cdot)$ is composed of only a few terms with respect to a large space of possible contributing terms. For instance, consider the non-linear function appearing in the Allen-Cahn equation 
\begin{equation}
    \mathbf{f}_p = p_1 \mathbf{u} + p_2\mathbf{u}^3 + p_3\mathbf{u}_{xx},
\end{equation}
where $\mathbf{p}=[5,\ -5,\ 0.0001]^\top$. Similarly, the non-linear function appearing in the Korteweg-de Vries equation is,
\begin{equation}
    \mathbf{f}_p = p_1 \mathbf{u} \mathbf{u}_x + p_2 \mathbf{u}_{xxx},
\end{equation}
where $\mathbf{p}=[-6,\ -1]^\top$.

The components of the parameter vector $\mathbf{p}$ can be computed via a least squares minimization formulation. This procedure often requires measuring $\mathbf{u}$ at $m$ different time points, and $n$ spatial locations. Consider that all these measurements are collected in a vector $\bU \in \R^{n\cdot m}$. Furthermore, if we assume that the function $\mathbf{f}_p(\cdot)$ have $k$ known linear, nonlinear, and partial derivatives terms, then it is possible to compute the following matrix:
\begin{equation}\label{dictionary}
	\Theta(\bU)= \big[\mathbf{1},\ \bU,\ \bU^2,\ \dots,\ \bU_x,\ \bU \bU_x,\ \dots \big], \ \ \Theta \in \R^{nm \times k},
\end{equation}
where each column of the matrix $\Theta$ contains all of the values of a particular term that constructs the right-hand side of \eqref{pde_1}, across all of the $n\cdot m$ space-time grid points where the data is collected. For example, if we consider measurements at $200$ spatial locations and $300$ time points with $\mathbf{f}$ comprising $5$ terms, then $\Theta \in \R^{200 \cdot 300 \times 5}$.

It is straightforward to compute the time derivative of  $\bU$ denoted by $\bU_t$, which is often implemented numerically. Having $\bU_t$ and other ingredients we can write the system of equation \eqref{pde_1} in the following form:
\begin{equation}\label{pde_2}
	\bU_t = \Theta(\bU) \mathbf{p},
\end{equation}
where $\mathbf{p} \in \R^{k}$ is an unknown parameter vector that has to be computed by a proper algorithm and its elements are coefficients corresponding to each term in the matrix $\Theta(\cdot)$ that describe the evolution of the dynamic system in time. The $i^{th}$ element of vector $\mathbf{p}$ is denoted by $p_i$.

The last step is to compute the parameter vector $\mathbf{p}$ through a least squares optimization procedure. A classical approach to solve \eqref{pde_2} comprise solving a system of linear equations known as the normal equations. This can be concretely expressed in the following fashion:
\begin{equation}\label{ls_sol}
\mathbf{p} = (\Theta^\top \Theta)^{-1} \Theta^\top \bU_t,
\end{equation}
which requires matrix inversion and multiplication of large size and it is not a good idea to explicitly form the inverse of a matrix, especially when dealing with large or ill-conditioned matrices \cite{ye2018accurate}. Therefore, in literature, authors often use regularized least squares (\texttt{RLS}) which is a family of methods for solving the least squares problem while using regularization to further constrain the resulting solution. Among \texttt{RLS} approaches, we can name \texttt{LASSO} \cite{friedman2001elements,tibshirani1996regression}, ridge regression \cite{kim2007interior}, and elastic net \cite{friedman2010regularization}. Note that in this article we make use of $\texttt{QR}$ factorization together with stochastic gradient descent technique to solve \eqref{pde_2}. Indeed, using the stochastic gradient descent algorithm which is a modification of the gradient descent method in the training loop of the \dnn, will alleviate the complexity of solution of least squares problem \eqref{ls_sol}, since we just use a random small part of our dataset instead of all of them. Hence, working on a portion of the dataset significantly reduces the computational load. We will elaborate on these settings in the upcoming sections.

\section{Greedy sampling method for \pde~parameter estimation}\label{sec:qdeim}
Estimating the parameters of \pde s is essential for understanding and predicting real-world phenomena, but it often poses significant challenges due to the inherent complexity and high dimensionality of these equations. One of these challenges is choosing the appropriate number of samples from a measured dataset which has an undisputed importance for the model recovery. 

In many real applications such as phase-field modeling or fluid dynamics, the state variables are often measured using sensors. The number of sensors is frequently constrained by physical or economic limitations, and the strategic placement of these sensors is crucial for achieving precise estimates. Regrettably, determining the ideal sensor locations to infer the \pde~parameters is inherently a combinatorial challenge, and existing approximation algorithms may not consistently produce effective solutions for all relevant cases. The topic of optimal sensor placement is a focus of research interest even in fields such as control theory and signal processing~\cite{both2021model}. Generally, five types of sensor placement approaches have been reported in the literature: (i) methods based on proper orthogonal decomposition (\propod) \cite{zhang2016pod} or compressed sensing \cite{moslemi2023sparse}, (ii) convex optimization methods \cite{joshi2008sensor}, (iii) greedy-based algorithms such as Frame-Sens \cite{ranieri2014near}, (iv) heuristic approaches such as population-based search \cite{lau2008tabu}, and (v) machine learning techniques \cite{wang2019reinforcement}.

Even though these approaches have shown good performance in tackling the curse of dimensionality by choosing the most informative $l$ locations from the assumed $n$ spatial ones, their applicability for the case of \pde~parameter estimation is limited due to lack of theoretical analysis (e.g. heuristic methods), lack of simplicity in implementations (machine learning approaches), and having conservative assumptions which might not be true in many cases (e.g. convex optimization methods).

More specifically, we are interested in feeding our most informative samples taken from the dataset into a \dnn~structure. 
It is worth highlighting that a \dnn~architecture has already shown its performance compared to traditional methods for the case of \pde~solution and \pde~discovery in \cite{both2021deepmod}. Unlike the traditional \pde~solvers that focus more on methods such as the \fdm~\cite{zeneli2021numerical} and \fem~\cite{bai2022local}, the \dnn~based approaches (such as \pinn) are meshfree and therefore highly flexible. Moreover, \dnn~has proven to regress along both the spatial and temporal axis for a given sample set \cite{both2021model}. However, as the the number of training samples increase for a DNN, the longer it takes to train the network. To alleviate this situation, it becomes crucial to choose a set of informative samples for training the network, and later, to estimate the parameters of the \pde.

In this article, we make use of \qdeim~as the sampling approach for a given dataset associated with the \pde s. In this regard, we give a brief introduction about \qdeim~by using its connection to a popular model order reduction approach named \propod. Consider the set of snapshots $\{u_1,\dots,u_m\}\in \R^n$ and an associated snapshot matrix $\cU=[u_1,\dots,u_m] \in \R^{n \times m}$ that is constructed by measuring the solution at $m$ different time points and $n$ different spatial locations of a \pde. In the conventional \propod,~we construct an orthogonal basis that can represent the dominant characteristics of the space of expected solutions that is defined as $\text{Range}\ \cU$, i.e., the span of the snapshots. We compute the singular value decomposition (\SVD) of the snapshot matrix $\cU$,
\begin{equation}\label{svd_1}
\cU = \bZ  \mathbf{\Sigma} \bY^\top,
\end{equation}

where $\bZ \in \R^{n\times n}$, $\mathbf{\Sigma} \in \R^{n \times m}$, and $\bY \in \R^{m \times m}$ with $\bZ^\top \bZ = \bI_n$, $\bY^\top \bY = \bI_m$, and $\mathbf{\Sigma}=\text{diag}(\sigma_1,\sigma_2,\cdots,\sigma_z)$ with $\sigma_1\ge \sigma_2 \ge \dots \ge \sigma_z \ge 0$ and $z = \min\{m,n\}$. The \propod~will select $\bV$ as the leading $\texttt{r}$ left singular vectors of $\bU$ corresponding to the $r$ largest singular values. Using \texttt{Python-Numpy} array notation, we denote this as $\bV = \bZ [:,:\texttt{r}]$. The basis selection via \propod~ minimizes $ \bV:= \min_{\Phi \in \R^{n\times r}} \| \cU - \Phi \Phi^\top \cU \|^2_F$, where $\|\cdot \|_F$ is the Frobenius norm, over all $\Phi \in \R^{n \times r}$ with orthonormal columns. In this regard, we can say $\cU = \bZ  \mathbf{\Sigma} \bY^\top \approx \bZ_{\texttt{r}} \mathbf{\Sigma}_{\texttt{r}} \bY_{\texttt{r}}^\top $, where matrices $\bZ_{\texttt{r}}$ and $\bY_r^\top$ contain the first $\texttt{r}$ columns of $\bZ$ and $\bY^\top$, and $\mathbf{\Sigma}_{\texttt{r}}$ contains the first $\texttt{r} \times \texttt{r}$ block of $\mathbf{\Sigma}$. More details regarding \propod~can be found in \cite{kunisch2002galerkin}.

Although the reduced-order model lies in the $\texttt{r}$-dimensional subspace, the conventional \propod~suffers the issue of lifting to the original space. Hence, in the literature, they tackle this issue with different approaches such as \deim~\cite{chaturantabut2010nonlinear}. An interesting advantage of \deim~is its flexibility to extend its results for the case of nonlinear function approximation beyond model order reduction. Moreover, the performance of the original \deim~algorithm has been improved by using \texttt{QR}-factorization in two aspects: (i) less upper bound error, (ii) simplicity and robustness in implementation.

\propod~has been used widely to select measurements in the state space that are informative for feature space reconstruction \cite{manohar2018data}. The method then has been called \qdeim. Indeed, we make use of the \qdeim~algorithm to select a set of most informative samples from a given snapshot matrix $\cU \in \R^{n \times m}$ of a \pde~for \pinn-based model discovery. \qdeim~takes advantage of the pivoted \texttt{QR} factorization and the \SVD, making it a powerful sampling approach. In particular, we consider \texttt{QR} factorization with column pivoting of $\bZ^\top_\texttt{r}$ and $\bY^\top_\texttt{r}$ to identify the most informative samples for the location and the time in the snapshot matrix $\cU$, respectively. The pivoting algorithm gives an approximate greedy solution approach for feature selection which is named submatrix volume maximization since matrix volume is defined as the absolute value of the determinant. \texttt{QR} column pivoting increases the volume of the submatrix constructed from the pivoted columns by choosing a new pivot column with maximal two-norm and then subtracting from every other column its orthogonal projection onto the pivot column. Note that \texttt{QR}-factorization has been implemented and optimized in most scientific computing packages and libraries, such as \texttt{MATLAB}, and \texttt{Python}. Further details about \qdeim, its theoretical analysis, and applications can be found in \cite{manohar2018data, chaturantabut2010nonlinear}.

\subsection{Applying \qdeim~algorithm on \pde~dataset}
We apply the \qdeim~ algorithm on the snapshot matrix $\cU$ to select the most informative samples in the spatiotemporal grid. To achieve this, we first perform a \texttt{SVD} of the snapshot matrix $\cU$ and compute matrices $\bZ$, $\mathbf{\Sigma}$, and $\bY^\top$. The selection of $\texttt{r}$ leading singular values can be done based on an appropriate precision value $\epsilon_{\texttt{thr}}$, which is related to the underlying dynamical system, and it can be chosen heuristically. $\epsilon_{\texttt{thr}}$ is often referred to as the energy criterion in literature \cite{gavish2014optimal}. In particular, the $\texttt{r}$ leading singular values are chosen such that the following quantity is satisfied:
\begin{equation}
1 - \frac{ \sum_{j=1}^{\texttt{r}} \sigma_j }{ \sum_{k=1}^{z} \sigma_k } < \epsilon_{\texttt{thr}} ,\ \  \texttt{r} < z,
\end{equation}  
Once the desired precision is achieved, we construct a reduced approximation of the snapshot matrix $\cU$ by using the first $\texttt{r}$ columns of matrix $\bZ$, the first $\texttt{r}$ singular values contained in the diagonal matrix $\mathbf{\Sigma}$, and the first $\texttt{r}$ rows of $\bY^\top$. Based on \texttt{Python} notation, we can represent this as $\bZ_\texttt{r} = \bZ[:,:\texttt{r}]$, $\mathbf{\Sigma}_\texttt{r} = \mathbf{\Sigma}[:\texttt{r},:\texttt{r}]$, and $\bY^\top_\texttt{r} = \bY^\top[:\texttt{r},:]$. To choose the important time and space indices, we apply \texttt{QR} decomposition with column pivoting on the reduced order matrices $\bY^\top_\texttt{r}$ and $\bZ^\top_\texttt{r}$, which contain the first $\texttt{r}$ left and right singular vectors. We represent the chosen indices corresponding to the most informative spatio-temporal points in the snapshot matrix $\cU$ by $\texttt{ind}_x$ and $\texttt{ind}_t$. To simplify the notation, we denote the pairs of space-time points by $(t_i,x_i)$, and the solution associated to it by $u_i$. \Cref{diemalgo} summarizes the core part of the \qdeim~sampling approach which takes the snapshot matrix $\cU$, spatio-temporal domains $x,\ t$ and the precision value $\epsilon_{\texttt{thr}}$ as the inputs and returns sampled pairs $(t_i, x_i)$ and corresponding measured value $u(t_i,x_i)$. We show by $\cN$ the cardinality of the sampled dataset.

\subsubsection{Exploiting locality in time}\label{time_division}
To capture the local dynamics of the dataset and select better points in the spatiotemporal domain, we decompose the time domain into equal intervals and apply \qdeim~on each sub-domain. We show the number of divisions for the time domain by $\texttt{t}_{\texttt{div}}$. In this setting, the subdomains do not overlap with each other, therefore the total number of selected points is the union of selected samples at each subdomain. For example, by using \texttt{Python} notation, if we divide the time domain into three parts, the first part can be written as $\cU [:,:\texttt{m}/3]$, the second part as $\cU [:,\texttt{m}/3:2\texttt{m}/3]$, and finally, the third part as $\cU [:,2\texttt{m}/3:]$. The reason behind this decomposition is to capture the local dynamics of the \pde~at each subdomain and sample the most informative portion of the snapshot matrix. More specifically, the system's behavior can vary across various domains, and certain physical attributes may exhibit notable distinctions. For instance, issues related to abrupt features such as shock waves may showcase these differences. On the other hand, by first dividing a sizable domain into smaller sub-domains, and then applying \qdeim~on each sub-domain independently, helps avoid the requirement of complex neural network structures for the \pde~parameter estimation.

\begin{algorithm}[t]
	\KwData{$\cU,\ \{x_k\}_{k=1}^{n},\ \{t_k\}_{k=1}^{m} ,\ \epsilon_{\texttt{thr}}$.}
	\KwResult{ $\texttt{ind}_t$, $\texttt{ind}_x$, domain sampled pairs $(t_i, x_i)$ and $u(t_i,x_i)$.
	}
	$\texttt{r} = 1 $\;
	
	$\bZ, \mathbf{\Sigma}, \bY^\top \gets$ \texttt{SVD}($\cU$), \Comment computing \texttt{SVD} on snapshot matrix $\cU$\;
	
	Find the lowest $ \texttt{r}$, such that $ 1 - \frac{\sum_{j=1}^{\texttt{r}} \sigma_j}{\sum_{j=1}^{z} \sigma_j}  \geq \epsilon_{\texttt{thr}}$;
	
	$\bZ_\texttt{r} \gets \bZ[:, : \texttt{r}]$, $\bY^\top_\texttt{r} \gets \bY^\top[:\texttt{r}, :]$, \Comment{selecting $\texttt{r}$ dominant left and right singular vectors}; 
	
	$\texttt{ind}_x \gets$ \texttt{QR}($\bZ^\top_\texttt{r}$, pivoting = \texttt{True}), \Comment{storing pivots from pivoted \texttt{QR} factorization of $\bZ^\top_\texttt{r}$};
	
	$\texttt{ind}_t \gets$ \texttt{QR}($\bY_\texttt{r}$, pivoting = \texttt{True}), \Comment{ storing pivots from pivoted \texttt{QR} factorization of $\bY^\top_\texttt{r}$};
	
	$x_i \gets \text{from}\ \texttt{ind}_x $, $t_i \gets \text{from}\ \texttt{ind}_t$, $u(t_i,x_i)$;
	
	\caption{Sample selection based on a two-way \qdeim~procedure}
	\label{diemalgo}
\end{algorithm}

\section{Greedy Sampling based Physics Informed Neural Network (\gspinn) for \pde~Parameter Estimation }%
\label{sec:method}
One notable constraint associated with \pinn s pertain to their extensive training costs, which can hinder their performance, particularly when addressing real-world applications that necessitate real-time execution of the \pinn~model. Consequently, it is vital to enhance the convergence speed of such models without compromising their efficacy. In this section, we propose our algorithm which blends \qdeim~as the sampling approach, \pinn~as the function approximation, and \texttt{QR}-factorization for the \pde~parameter estimation. In this article, we assume that we a priori know the terms that contribute to the \pde~dynamics. However, we do not have any information about their portions that are shown as coefficients, and our goal is to estimate them. For instance, consider the Allen-Chan equation with $\frac{\partial u}{\partial t} = p_1 u + p_2 u^3 + p_3 u_{xx}$ and associated unknown parameter vector $p^\top=[p_1,\ p_2,\ p_3]$.


\pinn s are primarily used for solving \pde s \cite{raissi2019physics, eshkofti2023gradient}. Unlike traditional numerical methods that discretize the domain and solve the equations on a grid, \pinn~treat the \pde~as a constraint in an optimization problem. The neural network is trained to minimize the discrepancy between the predicted and actual values of the \pde, effectively learning the underlying physics \cite{manavi2024trial}. An attractive feature of \pinn s is that it requires small modification for the case of \pde~parameter estimation \cite{lu2021deepxde}.

Another notable capability of the \pinn~is its Automatic Differentiation (AD). To elaborate more on this aspect, it is worth highlighting that numerical differentiation is generally performed by two primary approaches: finite differences or by using a spline interpolation method on the dataset followed by derivative calculation. Finite difference techniques work directly with the available data to compute the derivative. Although finite differences are computationally efficient and easily adaptable to higher dimensions, they are sensitive to noise and require closely spaced data points for accurate results. A more precise and commonly adopted alternative involves the application of spline interpolation for data differentiation. When employing splines for fitting, the data is approximated through a piecewise polynomial with ensured continuity at the boundaries. In practical terms, splines provide more accurate results. However, their scalability to higher dimensions, particularly when employing higher-order splines, is limited \cite{both2021model}. This limitation poses a challenge to model discovery, as it necessitates these higher orders to account for derivatives within the feature library. Using a fifth-order spline to approximate the data essentially translates to approximating the third-order derivative with only a second-order polynomial, thereby restricting its utility in the context of model discovery\cite{wandel2022spline}. While higher-order splines are limited to interpolation in a single dimension, \pinn~employs a neural network to interpolate along both the spatial and temporal axes. This liberates us from the requirement of using on-grid sampling, and we take advantage of this freedom by developing an alternative sampling technique. A combination of \deepnn~ with \sindy~based model discovery has been proposed in \cite{ both2021deepmod}, and several other works extended its results for the case of noisy and scarce dataset by using integration scheme in the training loop \cite{forootani2023, goyal2022neural}.

To setup our \pinn~formulation for the parameter estimation we first write the \pde~system of equations \eqref{pde_1} in a residual format:
\begin{equation}\label{residual_func}
\begin{aligned}
\mathcal{R}(u) & = - \frac{\partial u}{\partial t} - \mathbf{f}_p(\bu,\bu_x,\bu_{xx},\dots,x)\\
& = - \frac{\partial u}{\partial t} - \mathbf{\Theta}(u) p,
\end{aligned}
\end{equation}
and proceed to approximate $\bu(x,t)$ by a \dnn. In this regard we define a neural network $\cG_\theta(t,x)$ which has two inputs $x$ and $t$, and one output $\hat{\bu}$, i.e. $\hat{\bu} = \cG_\theta(t,x)$. Moreover, we harness the capabilities of \dnn, in particular, automatic differentiation (AD) to compute time derivative $\frac{\partial\hat{\bu}}{\partial t }$ with respect to its input variables, i.e. space $x$, and time $t$. Having these ingredients we define the following loss function:
\begin{equation}\label{objective_function}
\mathbf{\mathcal{L}} = \mu_1 \mathbf{\mathcal{L}}_{\texttt{MSE}} + \mu_2 \mathbf{\mathcal{L}}_{\texttt{deri}}, \ \ \ \mu_1,\mu_2 \in (0,1],
\end{equation}
where $\mathbf{\mathcal{L}}_{\texttt{MSE}}$ is the mean square error (\texttt{MSE}) of the output of the \dnn\ $\cG_\theta$ (denoted by $\mathbf{\hat{u}}$) with respect to its input sampled domain pairs $(t_i,x_i)$, and corresponding $u(t_i, x_i)$. The positive constants $\{\mu_{1}, \mu_{2} \}$ determine the weight of different losses in the total loss function. It is given as
\begin{equation}\label{mse_loss}
\mathbf{\mathcal{L}}_{\texttt{MSE}} = \frac{1}{\cN} \sum_{i=1}^{\cN}\Big\|{\bu(t_i,x_i) - \hat{\bu}(t_i,x_i)}\Big\|_2^2.
\end{equation}
$\mathcal{L}_{\texttt{MSE}}$ forces the \dnn\ to produce output in the vicinity of the measurements, and ${\mu_1}$ is its weight. $\mathbf{\mathcal{L}}_{\texttt{deri}}$ is the residual \eqref{residual_func} and aims to compute the parameter vector $p$. It is computed as follows:
\begin{equation}\label{reg_loss}
\mathbf{\mathcal{L}}_{\texttt{deri}} = \frac{1}{\cN} \sum_{i=1}^{\cN} \Big\| \frac{\partial \hat{\bu}(t_i,x_i)}{\partial t_i} - \mathbf{\Theta}\big( \hat{\bu}(t_i,x_i) \big) \hat{p} \Big\|_2^2,
\end{equation} 
where $\hat{p}$ is replaced by computing the $\texttt{QR}$-factorization of the matrix $\mathbf{\Theta}$ and its pseudo-inverse as follows 
\begin{equation}\label{param_vec}
\begin{aligned}
\bQ \bR & = \mathbf{\Theta}(\hat{\bU}), \\
\hat{p} & = \bR^{-1} \bQ^\top \hat{\bU}_t,
\end{aligned}
\end{equation}
here $\hat{\bU}$ and $\hat{\bU}_t$ are column vectors with components $ \hat{\bu}(t_i,x_i)$ and $ \frac{\partial\hat{\bu}(t_i,x_i)}{\partial t_i}$, respectively. It is worth highlighting that the parameter vector $p$ will be updated alongside the weights and biases of the \dnn, and the matrix $\Theta$ is calculated by \eqref{dictionary}. \Cref{alg_two} summarizes the procedure that has been explained in this section.

\alglanguage{pseudocode}
\begin{algorithm}[h]
\KwData{$\cU$, $x$, $t$,  $\epsilon_{\texttt{thr}}$ for the \qdeim~ algorithm, matrix $\mathbf{\Theta}$, a neural network $\cG_\theta$ (parameterized by $\theta$), maximum iterations ${\texttt{max-iter}}$, and parameters $\{\mu_1,\mu_2\}$.}
\KwResult{ Estimated parameter vector $p$, defining governing equations}

$(t_i, x_i),\ u(t_i,x_i) \gets$ Apply \qdeim($\cU$) based on \Cref{diemalgo} \Comment{selecting most informative samples} \;


Initialize the \dnn~module parameters, and the parameter vector $p$\;
$k=1$\;

\While{ $k<$ \texttt{max-iter}  }{
Feed the domain pairs $(t_i,x_i)$ as an input to the \dnn~($\cG_\theta$) and predict output $\mathbf{\hat{u}}(t_i,x_i)$ \;

Compute the derivative information $\frac{\partial \mathbf{\hat{u}}(t_i,x_i)}{ \partial t_i }$ using automatic differentiation \;
	
Compute the cost function~\eqref{objective_function} \;

Update the parameters of \dnn\ ($\theta$) and the parameter vector $p$ using gradient descent \;
}

\caption{Greedy Sampling based Physics Informed Neural Network (\texttt{GS-PINN}) for \pde~Parameter Estimation}
\label{alg_two}
\end{algorithm}

\begin{figure}[t]
	\centering
	\includegraphics[width=\textwidth]{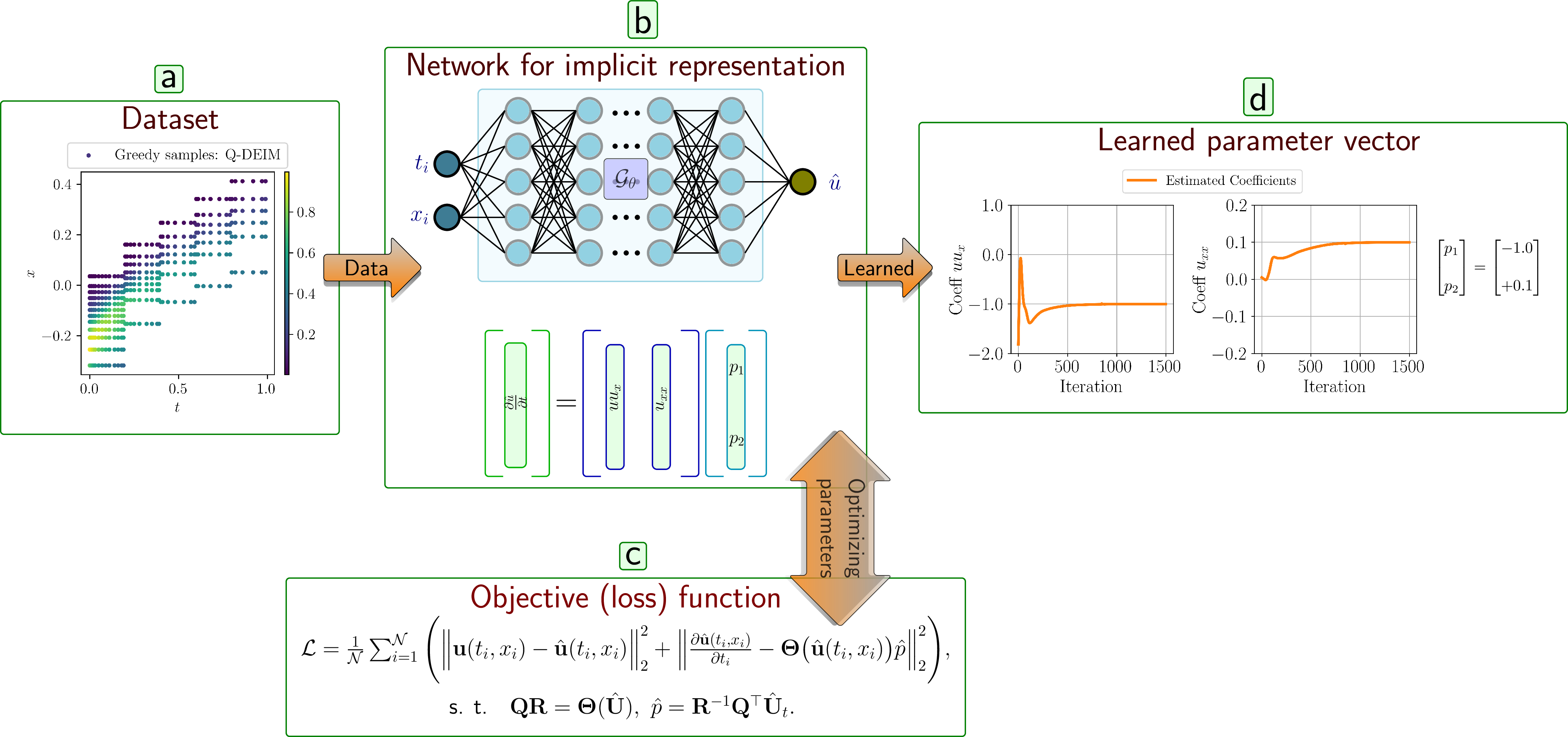}
	\caption{A schematic diagram of \gspinn~for \pde~parameter estimation. (a) sampling \pde~dataset with \qdeim~algorithm, (b) feeding the pairs $(t_i,x_i)$ time and space resulted from \qdeim~into \dnn, (c) using the output of the \dnn~as the function approximator, construct matrix $\cU$, (d) estimating the parameters of the \dnn~with stochastic gradient descent and parameter vector $p$ by considering the loss function $\cL$.
	}
	\label{iNeuralSINDy_schematic}	
\end{figure}


\section{Simulation Results} \label{sec:Simulation}
In this section we provide several simulation examples for the case of \pde~parameter estimation with our proposed \gspinn~ \pde~estimator. We also provide a comparison for the case where the spatio-temporal dataset is randomly sampled. The simulation examples have different level of complexity and non-linearity. To have a quantitative comparison and evaluate the performance of \gspinn~versus random sampling, we consider the following relative error for each element of parameter vector $p_i$:
\begin{equation}\label{error_criteria}
\texttt{error}_{{p_i}}=\frac{ \Big|p^{\texttt{truth}}_i - p^{\texttt{est}}_i  \Big| }{ \Big| p^{\texttt{truth}}_i \Big |}
\end{equation}
where $| \cdot |$ denotes the absolute value; $p^{\texttt{truth}}_i$ and $p^{\texttt{est}}_i$ are the true model coefficient and the estimated coefficient, corresponding to the coefficient ${p_i}$, respectively. The main point to define such criteria is the fact that in some \pde s the scale of coefficients for different terms are not in the same range (such as Allen-Cahn equation with $\frac{\partial u}{\partial t}=0.0001 u_{xx} - 5u^3 + 5u$), therefore it is logical to quantify the errors for each element of the parameter vector separately. 

In this article we examine the parameter estimations for three \pde s: (i) the Allen-Cahn equation, (ii) the Burgers' equation, and (iii) the Korteweg-de Vries (KdV) equation. We recall that in our setting we do not consider the boundary conditions of \pde s, since our goal is not solving them, rather we desire to estimate their parameters. For each example, the simulation is divided into two parts. In the first part, we employ \qdeim~on the given snapshot matrix for a \pde~to identify the most informative samples. This step requires to have a trade off between the number of divisions in time domain $\texttt{t}_{\texttt{div}}$ and precision value $\epsilon_{\texttt{thr}}$, since we can manipulate the cardinality of our sample set. Then we apply \Cref{alg_two} to estimate the parameters of the \pde. In the second part, we do a comparison between \qdeim~sampling and random sampling approach based on the same sample size. To identify the minimum and maximum number of samples for the comparison, in case of the \qdeim~algorithm, the time domain division is considered to have integer values $\texttt{t}_{\texttt{div}}=1,\dots,4$ and we choose the precision value $\epsilon_{\texttt{thr}}$ to be equally spaced logarithmic values, spanning the range from a minimum value to maximum value, e.g., $20$ equally spaced logarithmic values, spanning the range from $10^{-10}$ to $10^{-2}$. This approach is commonly used in scientific and computational applications, such as when specifying tolerance levels for numerical calculations or for testing a range of values in a logarithmic scale. Having the collection of samples and for each pair $(\texttt{t}_{\texttt{div}},\ \epsilon_{\texttt{thr}})$ we apply our \texttt{GS-PINN} \Cref{alg_two} for \pde~parameter estimation. By doing so it reveals the effect of time domain division $\texttt{t}_{\texttt{div}}$, the precision value $\epsilon_{\texttt{thr}}$ and the sample size on the \texttt{GS-PINN} performance for parameter estimation. To choose the sample size for the simulation of random sampling approach we make use of sample range computed via \qdeim~algorithm for pairs $(\texttt{t}_{\texttt{div}},\ \epsilon_{\texttt{thr}})$ as explained earlier for each \pde. In the sense that we select the minimum and the maximum value among all sample sizes that we identified by \qdeim~algorithm for pairs $(\texttt{t}_{\texttt{div}},\ \epsilon_{\texttt{thr}})$ and we start from the minimum sample size and we reach to the maximum value with $10$ linearly spaced steps. For each sample size we pick samples from \pde~snapshot matrix randomly and we employ our \texttt{PINN} setup. We plan to conduct simulations $5$ times for each sample size without specifying a random seed, resulting in a total of $55$ experiments. For the sake of simplicity the \pinn~architecture which make use of random sampling is named Random Sampling \pinn~(\rspinn). The details of each simulation example will be provided accordingly.

\paragraph{Data generation.}
We acquired our dataset associated to each \pde~from the repository provided by the works reported in \cite{rudy2017data, both2021deepmod}. Moreover, we perform a data-processing step before feeding them to the \qdeim~and \dnn~ for both time $t$ and the space $x$ domains. In particular, we map the time domain $t$ into the interval $[0,\ 1]$ and the space domain $x$ into the interval $[-1,\ 1]$. For each \pde~we also mention the original range of the space-time domain dataset.

\paragraph{Architecture.}
We utilize multi-layer perceptron networks with periodic activation functions, particularly embracing the \texttt{SIREN} architecture as introduced by \cite{sitzmann2019siren}. This approach allows us to obtain an implicit representation from measurement data, and we will customize the number of hidden layers and neurons for each specific example. Our \dnn~ architecture for all the examples have $3$ hidden layers, each having $128$ neurons.

\paragraph{Hardware.}
In our neural network training and parameter estimation efforts to uncover governing equations, we employed an \nvidia RTX A4000 GPU with 16 GB RAM. For CPU-intensive tasks such as data generation, we harnessed the power of a 12th Gen \intel~\coreifive-12600K processor equipped with 32 GB RAM.

\paragraph{Training set-up.}
We use the Adam optimizer with $\texttt{learning\_rate}=10^{-5}$ to update the parameter vector $p$ that is trained alongside the \dnn\ parameters \cite{kingma2014adam}. In particular we employ \texttt{CyclicLR}, which is a learning rate scheduling technique commonly use in deep learning with the \texttt{PyTorch} framework \cite{paszke2019pytorch}. It allows for dynamic adjustments of the learning rate during training to potentially improve the training process. The key parameters include $\texttt{base}\_{\texttt{lr}}= 0.1\times \texttt{learning\_rate}$ and $\texttt{max\_lr}= 10\times \texttt{learning\_rate}$, which define the lower and upper bounds of the learning rate respectively, $\texttt{cycle}\_{\texttt{momentum}}$ to control cycling of momentum, and \texttt{mode} to specify the learning rate cycling strategy. In our setting, it employs the $\texttt{exp}\_{\texttt{range}}$ mode, which means the learning rate oscillates exponentially between the specified lower and upper bounds over a set number of iterations as determined by $\texttt{step\_size\_up}=1000$. This technique can be effective in training deep neural networks by promoting faster convergence and potentially helping to escape local minima during optimization. The positive constant parameters $\mu_{1,2}$ are considered to have equal value, i.e. $\mu_{1,2}=1$. Number of time domain division ($\texttt{t}_{\texttt{div}}$) to apply \qdeim~for each \pde~dataset, precision value ($\epsilon_{\texttt{thr}}$), and maximum number of iterations (${\texttt{max-iter}}$) will be mentioned separately for each example.


\subsection{Allen-Cahn equation}%
\label{subsec:AC}
The Allen-Cahn equation lies at the heart of modelling the transformation of a physical quantity, often called the order parameter, as a material undergoes a transition from one phase to another. This equation elegantly encapsulates phenomena such as solidification, crystal growth and the emergence of domain patterns in magnetic materials. It belongs to the broader class of reaction-diffusion equations, which relate diffusion processes to a double-well potential energy function. Consequently, the equation encapsulates the inherent drive of the system to minimise its free energy, leading to the formation of domains with distinctive properties. These domains are characterised by smooth variations in the order parameter, with sharp transitions occurring at domain walls. Their versatility extends to a wide range of applications, including the study of grain boundaries in materials science, the modelling of phase separation in binary fluids, and the elucidation of the intricacies of pattern formation in neural networks in neuroscience. The paramount importance of understanding phase transitions and the intricate development of complex patterns in complex patterns in physical systems underlines the invaluable role played by the Allen-Cahn equation \cite{feng2003numerical}.

Let us extract the Allen-Cahn equation by using Ginzburg-Landau free energy \cite{wight2020solving}
\begin{equation}
\cF=  \int_{w}^{} \frac{\gamma_1}{2} \big| \grad u \big|^2 + \frac{\gamma_2}{4} (u^2 -1) \,dx ,
\end{equation}
where $\gamma_1$ and $\gamma_2$ are parameters. By computing $L^2$ gradient flow, we obtain the Allen-Cahn equation
\begin{equation}\label{AC_eq}
\begin{aligned}
&\frac{\partial u}{\partial t} + \gamma_1\ \Delta u + \gamma_2\ (u-u^3)=0, \\
&\frac{\partial u}{\partial t} = -\gamma_1\ u_{xx} - \gamma_2\ (u-u^3).
\end{aligned}
\end{equation}

\begin{figure}[!ht]
	\centering
	\includegraphics[width=1\textwidth]{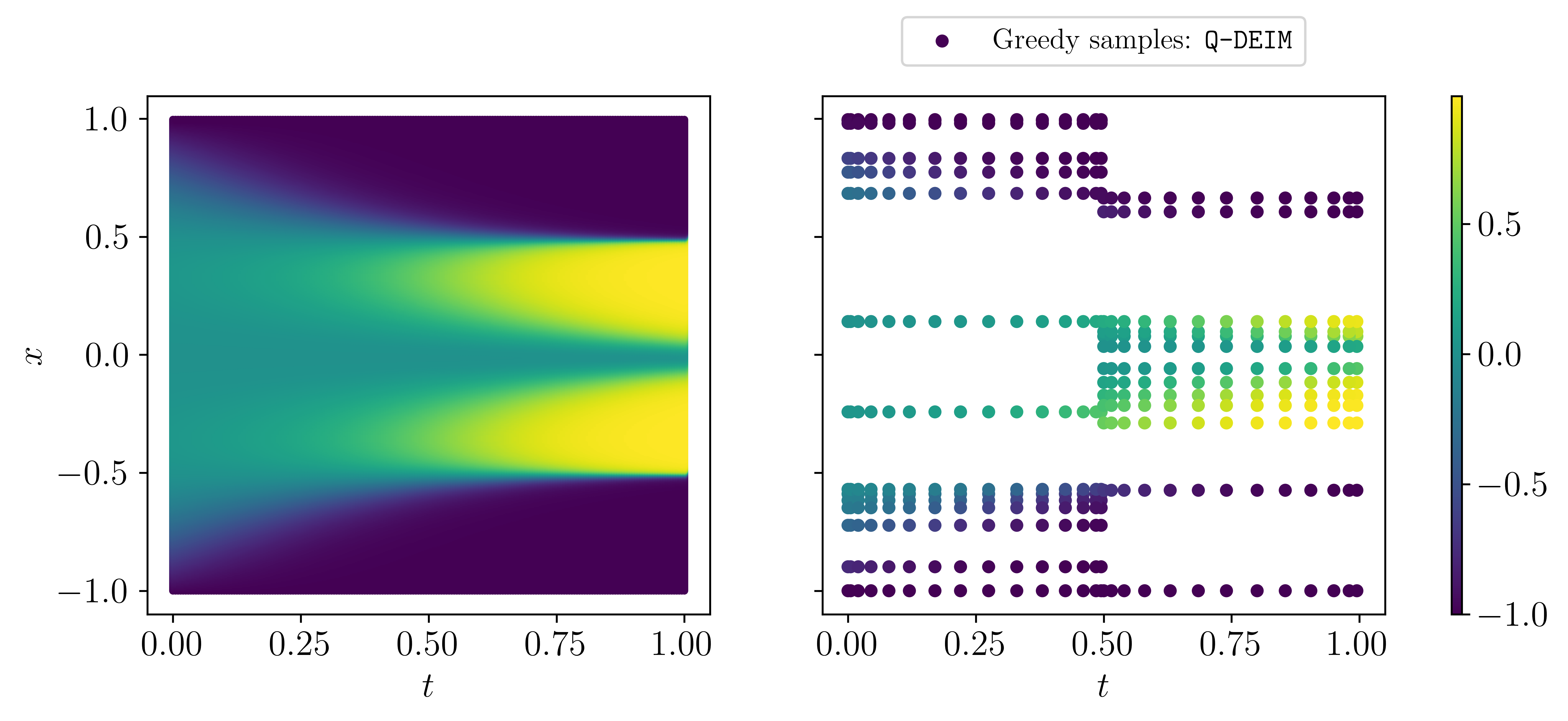}
	\caption{(left) Entire dataset; (right) Greedy samples by \qdeim~algorithm for Allen-Cahn equation  (\Cref{subsec:AC})}
	\label{AC_qdeim}	
\end{figure}

\begin{figure}[!h]
	\centering
	\includegraphics[width=0.8\textwidth]{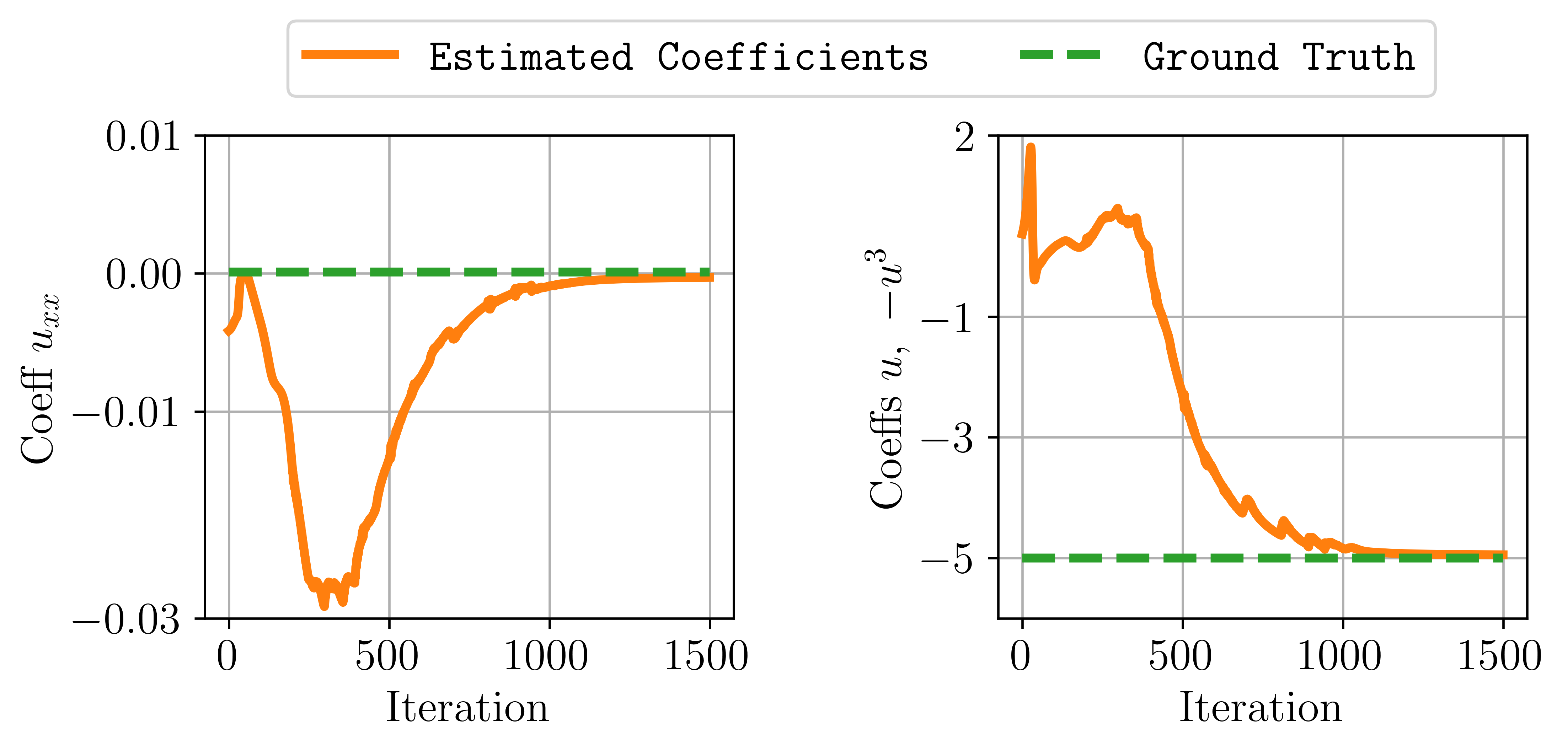}
	\caption{ Estimated coefficients by \texttt{GS-PINN} for Allen-Cahn equation (\Cref{subsec:AC})}
	\label{AC_qdeim_coeffs}	
\end{figure}

We employed \texttt{GS-PINN} for \pde~parameter estimator on a dataset corresponding to Allen-Cahn equation with periodic boundary conditions in one dimension with nominal parameters $\gamma_1=0.0001$, and $\gamma_2=5$. In this example the original dataset is taken from the work reported in \cite{wight2020solving}. The space domain $x$ and time domain $t$ were already mapped into the interval $x \in [-1,\ 1]$ and $t \in[0,\ 1]$ and were discretized into $n=512$ locations and $m=201$ points, respectively. Therefore the snapshot matrix $\cU \in \mathbb{R}^{n\times m}$ has $102912$ elements which represents the curse of dimensionality of training the \pinn~architecture for the \pde~parameter estimation. The snapshot matrix of the original dataset is shown in the \Cref{AC_qdeim} (left) which represents the process of phase separation in time. To apply \qdeim~algorithm we set the time domain division $\texttt{t}_{\texttt{div}}=2$, precision value $\epsilon_{\texttt{thr}}= 10^{-8}$. In the \Cref{AC_qdeim} (right) we see the most informative samples that are selected with \qdeim~algorithm. In total, $394$ samples are selected to be utilized in \pinn~parameter estimation with maximum number of iterations ${\texttt{max-iter}}=1500$. The trajectory of the coefficients estimated by \gspinn~are depicted in \Cref{AC_qdeim_coeffs}. The parameter vector is $p=[4.95, -4.95, -0.0003]^\top$ which based on \eqref{AC_eq} $\gamma_1=-p_3$ and $\gamma_2=-p_1=-p_2$. We observe that with having $0.383 \%$ of the entire dataset (i.e. $\frac{394}{102912}\approx 0.00383\times 100= 0.383\%$) we earned acceptable precision. We notice that \gspinn~can not estimate the parameter $\gamma_1$ properly since its value is considerably less than $\gamma_2$ in the nominal equation.

\begin{figure}[!h]
	\centering
	\includegraphics[width=0.8\textwidth]{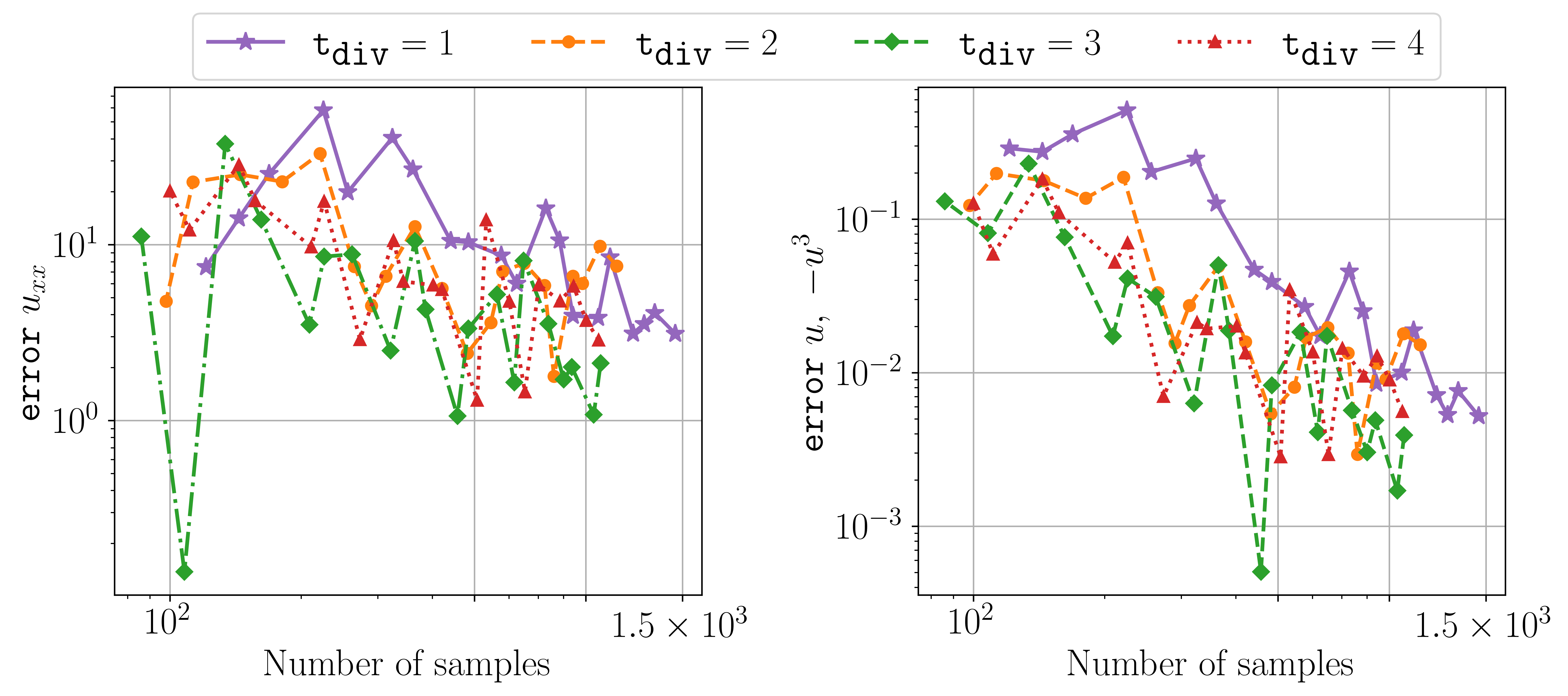}
	
	\caption{Performance of \gspinn~to estimate coefficients corresponding to Allen-Cahn \pde~equation (\Cref{subsec:AC}) with different time divisions $\texttt{t}_{\texttt{div}}$ and precision values $\epsilon_{\texttt{thr}}$}
	\label{AC_qdeim_sensitivity_1}	
\end{figure}

In the second scenario, i.e. comparing \gspinn~with \rspinn, we choose $20$ equally spaced logarithmic values for precision $\epsilon_{\texttt{thr}}$, spanning the range from $10^{-13}$ to $10^{-4}$ corresponding to each time division $\texttt{t}_{\texttt{div}}=1,\ 2,\ 3,\ 4$. The minimum number of samples corresponds to the pair $(\texttt{t}_{\texttt{div}},\ \epsilon_{\texttt{thr}})= (3,\ 10^{-4})$ with $86$ samples, while the maximum number of samples corresponds to $(\texttt{t}_{\texttt{div}},\ \epsilon_{\texttt{thr}})= (1,\ 10^{-13})$ with $1444$ samples. The performance of \gspinn~for each pair $( \texttt{t}_{\texttt{div}},\ \epsilon_{\texttt{thr}})$ is shown in the figure \Cref{AC_qdeim_sensitivity_1}. The first notable problem that we recognize is that the \gspinn~does not have a good estimation about coefficient corresponding to the term $ u_{xx}$, however \gspinn~is able to recover the coefficients corresponding to the terms $u$ as well as $u^3$. As we can see $ \texttt{t}_{\texttt{div}}=1$ has the worst performance with both highest relative error and highest sample range. In the sense that for the fixed relative error  $\texttt{t}_{\texttt{div}}=1$ almost requires more samples compare to the $ \texttt{t}_{\texttt{div}}=2,\ 3, \ 4$. This result reveals the necessity of finding the suitable value for $\texttt{t}_{\texttt{div}}$ to have a less relative error in \pde~parameter estimation with \gspinn~and therefore better performance. Moreover, we see that for the case of $\texttt{t}_{\texttt{div}}=3$ and $\texttt{t}_{\texttt{div}}=4$ the sample range and error values are almost the same and it is hard to distinguish which one performs better. This shows that with a fixed precision value $\epsilon_{\texttt{thr}}$ increasing the time division $\texttt{t}_{\texttt{div}}$ has no sensible effect on the quality of the estimation for Allen-Cahn \pde~equation. It is worth to highlight that the minimum relative error corresponds to the pair $(\texttt{t}_{\texttt{div}}, \epsilon_{\texttt{thr}})=(3,\ 1.833\times 10^{-9})$ with $457$ samples.

\begin{figure}[!h]
	\centering
	\includegraphics[width=0.8\textwidth]{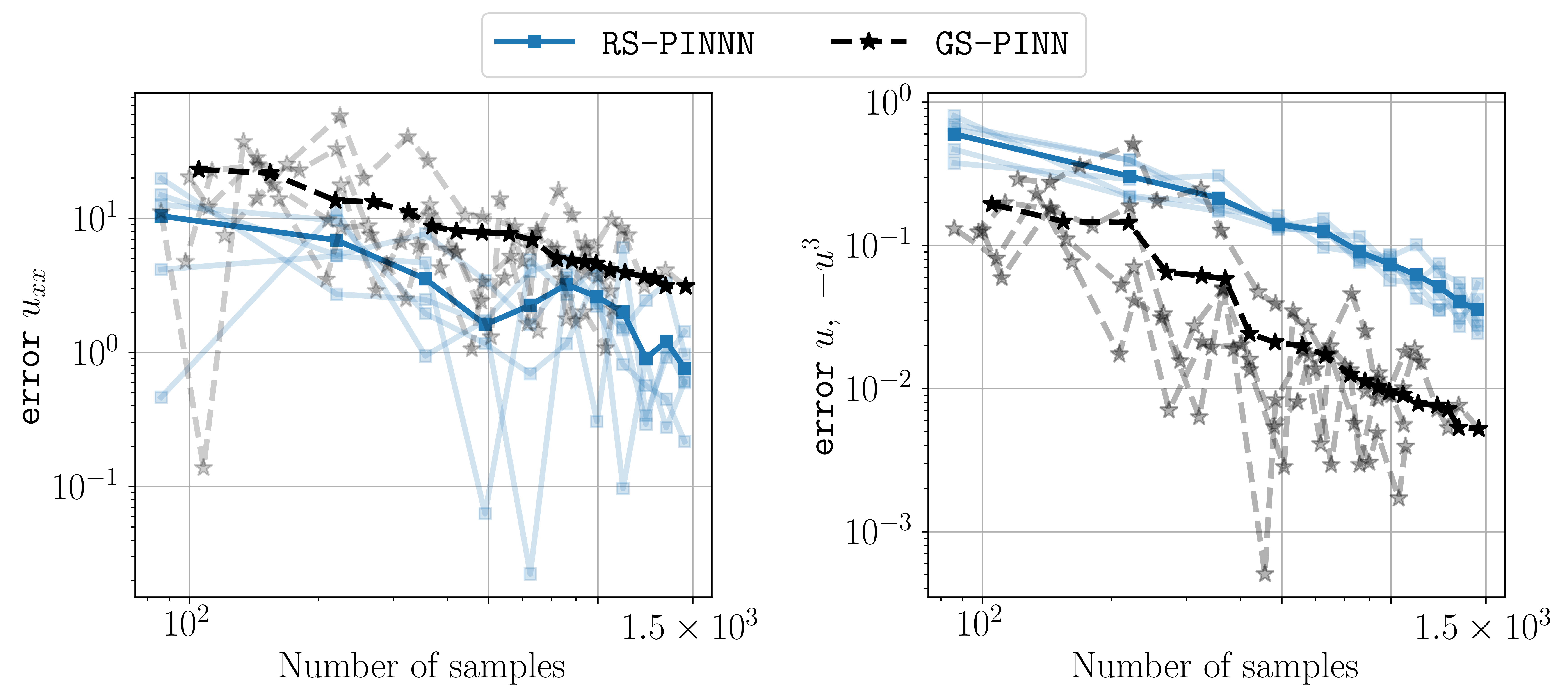}
	\caption{Comparing \gspinn~with \rspinn~to estimate coefficients corresponding to Allen-Cahn \pde~equation (\Cref{subsec:AC}) with different time divisions $\texttt{t}_{\texttt{div}}$ and precision values $\epsilon_{\texttt{thr}}$}
	\label{AC_qdeim_sensitivity_2}	
\end{figure}

Having the sample range computed by \gspinn~we are ready to evaluate the performance of \rspinn~to recover the \pde~coefficients. The result of this comparison between \gspinn~and \rspinn~is shown in the \Cref{AC_qdeim_sensitivity_2}. For the sake of simplicity in the plots related to the performance of \rspinn, we take the average value of the relative errors for a given fixed number of samples associated to $5$ experiments. Since the range of samples is related to the time division $\texttt{t}_{\texttt{div}}$, for \gspinn~we employ \texttt{K-means} clustering with $k=20$, see \Cref{sec:appendix} and \cite{arthur2007k}, to find the nearest cluster centroid corresponding to each pair sample size and relative error. This significantly simplifies the evaluation of \gspinn~and \rspinn~under different configurations. These curves are shown in the \Cref{AC_qdeim_sensitivity_2} for both \gspinn~(dashed line) and \rspinn~(solid line). We see that \rspinn~can not estimate the coefficient corresponding to the term $u_{xx}$, however it has lower relative error with respect to \gspinn. This is due to the fact that the \dnn~is not capable to reach a small precision value in the range $0.0001$ which coefficient corresponding to the term $u_{xx}$ is. Regarding the estimated values for the terms $u$ and $u^3$, we see that \gspinn~outperforms \rspinn~in all the settings except in one case.

\subsection{Burgers' equation}
\label{subsec:Burger}
Burgers' equation captures the dynamics of a quantity, typically denoted as the velocity or density, within a medium characterized by convection and diffusion processes. It represents a crucial intermediary between the simpler linear convection-diffusion equation and the more complex nonlinear conservation equations. The equation is characterized by its capacity to describe a multitude of phenomena, including the formation of shocks, rarefaction waves, and the emergence of complex wave patterns. Within the realm of fluid dynamics, Burgers' equation finds application in modeling traffic flow, the dynamics of shock waves in compressible fluids, and the behavior of dispersive waves in shallow water. Furthermore, it has significance in other fields such as plasma physics, acoustics, and the study of nonlinear optical systems \cite{ozics2003finite}. The one-dimensional form of the Burgers' Equation is typically written as:
\begin{equation}\label{Burger}
\frac{\partial u}{\partial t} + \lambda\ u u_x + \nu\ u_{xx}=0.
\end{equation}

In order to study the Burgers' equation, as described in \cite{rudy2017data}, it is discretized over the domain \(x \in [-8, 8]\) into $n=256$ points and in a time interval \(t \in [0, 10]\) with $m = 101$ time steps. 

\begin{figure}[!tb]
	\centering
	\includegraphics[width=1\textwidth]{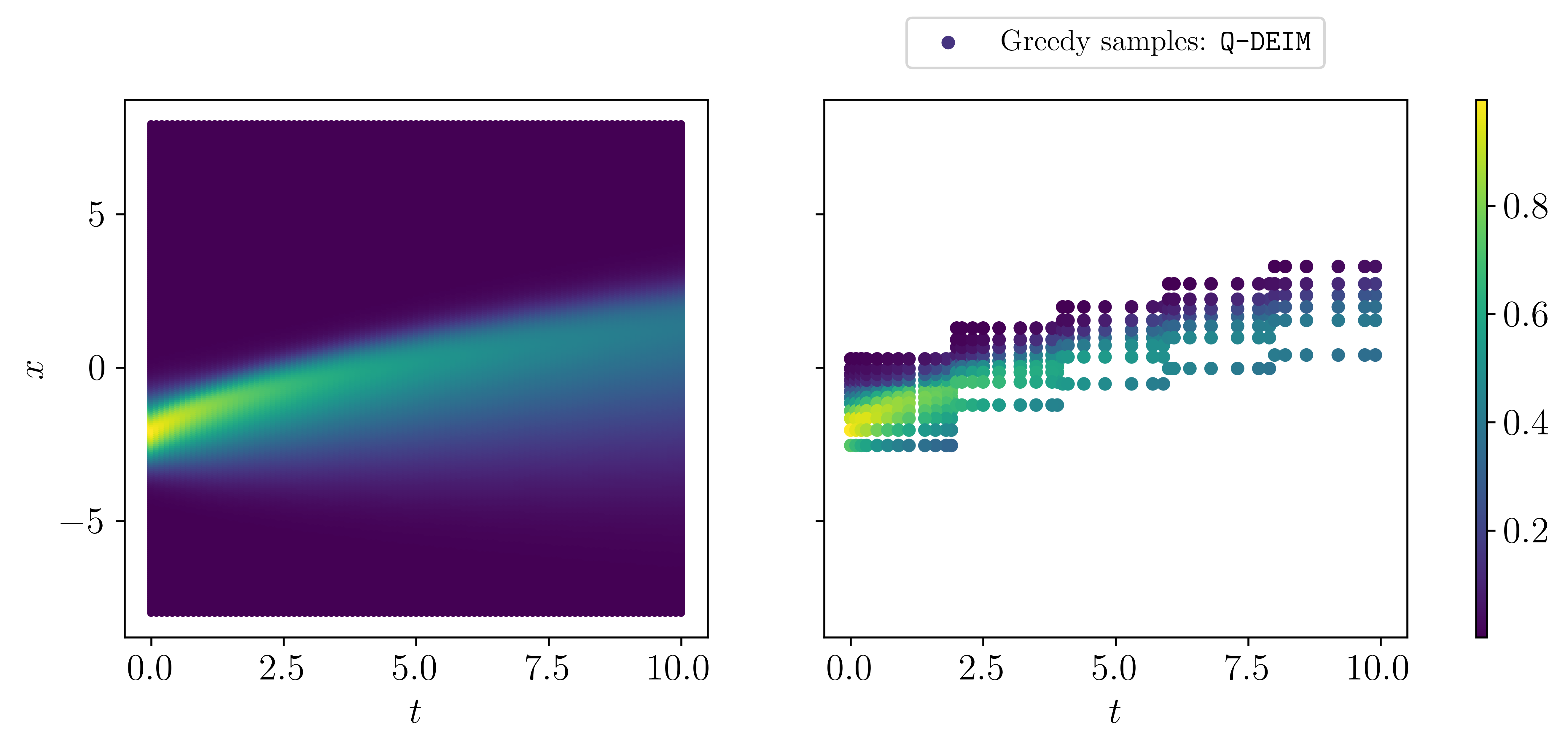}
	\caption{(left) Entire dataset ; (right) Greedy samples by \qdeim~algorithm for Burger's equation  (\Cref{subsec:Burger})}
	\label{Burger_qdeim}	
\end{figure}

Hence, the snapshot matrix $\cU$, with dimensions $n\times m$, comprises a total of $25,856$ elements, highlighting the challenge posed by the curse of dimensionality when training the \pinn~architecture for \pde~parameter estimation. The left side of \Cref{Burger_qdeim} displays the snapshot matrix of the original dataset, illustrating the Burgers' equation dynamic. To employ the \qdeim~algorithm, we specify a time domain division of $\texttt{t}_{\texttt{div}}=5$ and a precision threshold value of $\epsilon_{\texttt{thr}}= 10^{-6}$.

\begin{figure}[!tb]
	\centering
	\includegraphics[width=0.8\textwidth]{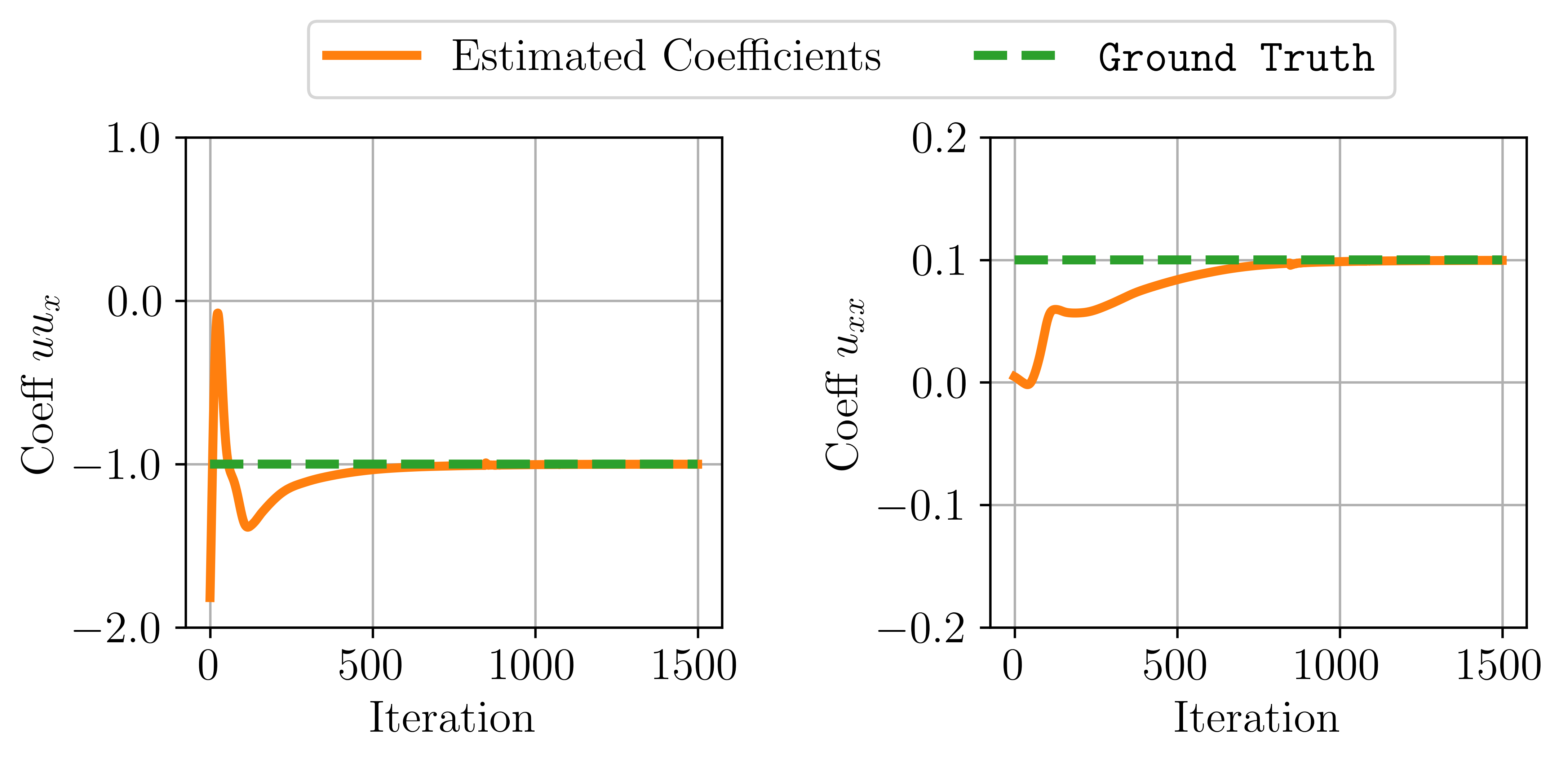}
	\caption{ Estimated coefficients by \texttt{GS-PINN} for Burger's equation (\Cref{subsec:Burger})}
	\label{Burger_qdeim_coeffs}	
\end{figure}

In \Cref{Burger_qdeim} on the right-hand side, we observe the selection of the most informative samples accomplished using the \qdeim~algorithm. A total of $359$ samples have been chosen for utilization in \gspinn~parameter estimation, with a maximum iteration set at ${\texttt{max-iter}}=1500$. The evolution of the coefficients estimated by the \gspinn~method is illustrated in \Cref{Burger_qdeim_coeffs}. The computed parameter vector $p$ after $1500$ iteration is $p^\top = [-1.000,\   0.0996]$ where based on \eqref{Burger} we have $\lambda = -p_1$ and $\nu = -p_2 $. We see that with approximately $1.39\%$ of the entire dataset we can recover Burgers' \pde~parameters.

\begin{figure}[!tb]
	\centering
	\includegraphics[width=0.8\textwidth]{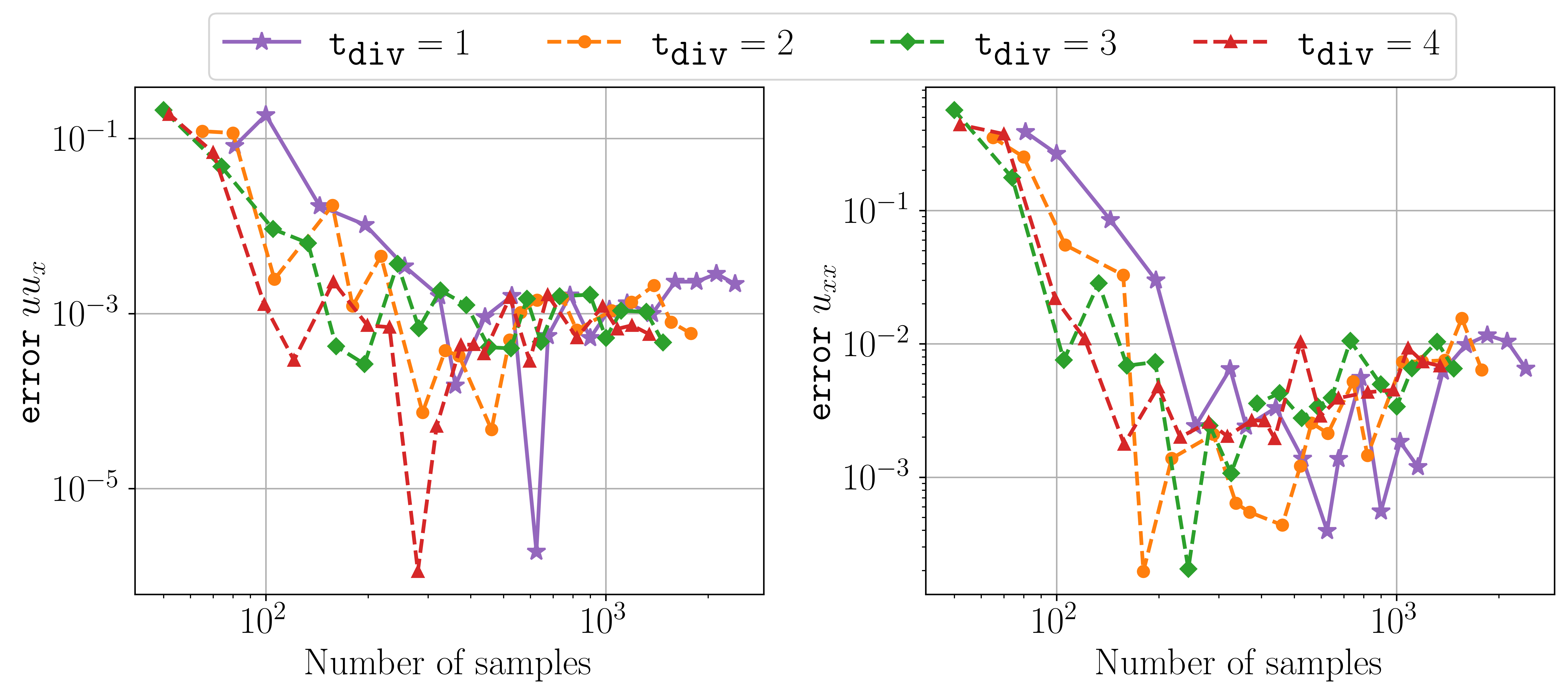}
	\caption{ Performance of \gspinn~to estimate coefficients corresponding to Burger's \pde~equation (\Cref{subsec:Burger}) with different time divisions $\texttt{t}_{\texttt{div}}$ and precision values $\epsilon_{\texttt{thr}}$}
	\label{Burger_qdeim_sensitivity_1}	
\end{figure}

In \Cref{Burger_qdeim_sensitivity_1} the relative error corresponding to different pairs $( \texttt{t}_{\texttt{div}}, \epsilon_{\texttt{thr}} )$ is shown. We select $20$ equally spaced logarithmic values for precision $\epsilon_{\texttt{thr}}$, spanning the range from $10^{-10}$ to $10^{-2}$ corresponding to each time division $\texttt{t}_{\texttt{div}}=1,\ 2,\ 3,\ 4$. The minimum number of samples corresponds to the case of $\texttt{t}_{\texttt{div}}=3$ with associated sample size $50$. We see that $\texttt{t}_{\texttt{div}}=1$ has the worst performance in the interval with less number of samples, however it improves with more number of samples and in some cases outperforms the others. Moreover, in the interval corresponding to higher sample size manipulating $\texttt{t}_{\texttt{div}}$ has no sensible effect on the performance of \gspinn~; this is more evident for the relative error corresponding to the term $uu_x$ in \Cref{Burger_qdeim_sensitivity_1}. The lowest relative error associated to the term $u_{xx}$ corresponds to the pair $(\texttt{t}_{\texttt{div}},\ \epsilon_{\texttt{thr}})=(2,\ 2\times 10^{-4})$, while for the term $uu_x$ lowest relative error corresponds to the pair $(\texttt{t}_{\texttt{div}},\ \epsilon_{\texttt{thr}})=(4,\ 1.13 \times 10^{-05})$.

\begin{figure}[!tb]
	\centering
	\includegraphics[width=0.8\textwidth]{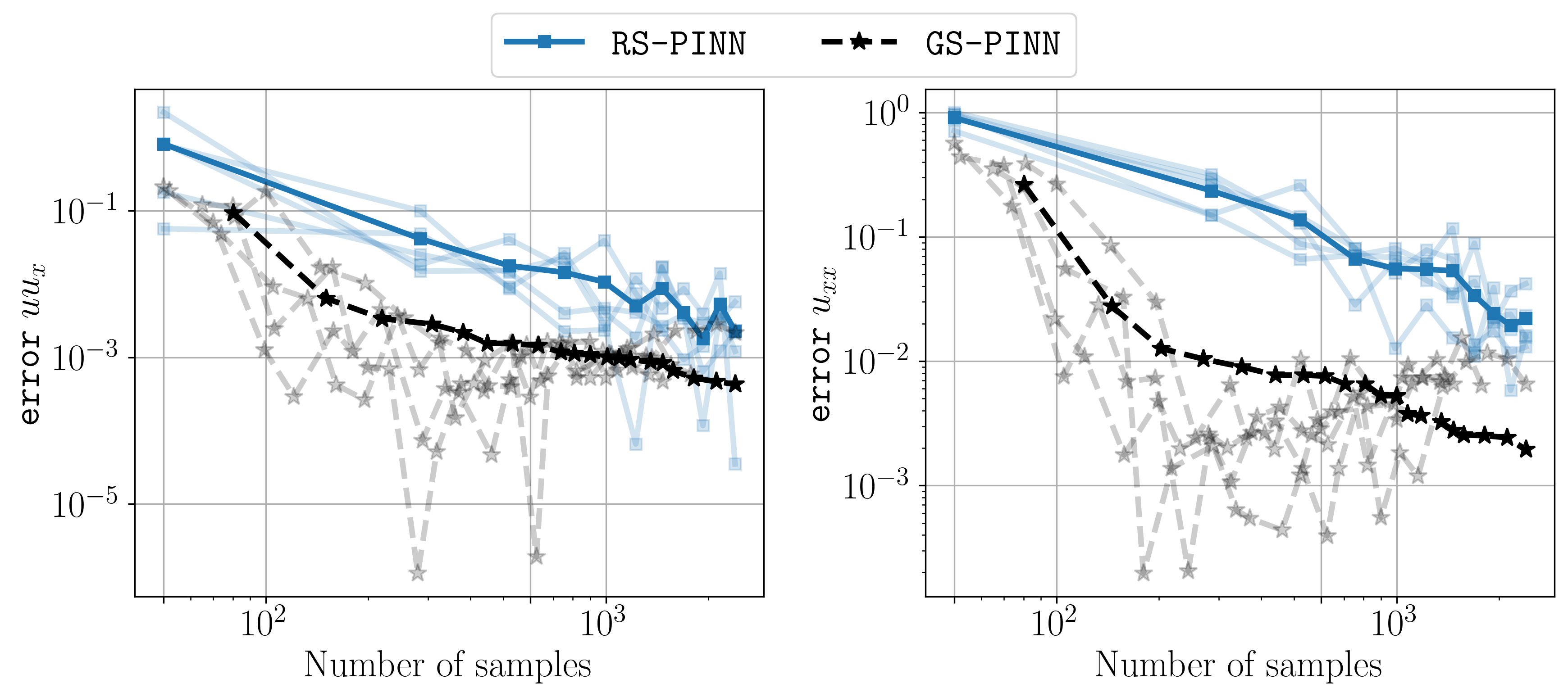}
	\caption{Comparing \gspinn~with \rspinn~to estimate coefficients corresponding to Burger's \pde~equation (\Cref{subsec:Burger}) with different time divisions $\texttt{t}_{\texttt{div}}$ and precision values $\epsilon_{\texttt{thr}}$}
	\label{Burger_qdeim_sensitivity_2}	
\end{figure}

With the sample range calculated using \gspinn, we are now prepared to assess how well \rspinn~performs in the recovery of \pde~coefficients. The outcome of this comparison between \gspinn~and \rspinn~can be seen in \Cref{Burger_qdeim_sensitivity_2}. To simplify the plots related to the performance of \rspinn, we take the average value of the relative errors across five experiments for a fixed number of samples.

Like the previous example, as the sample range is linked to the time division parameter, denoted as $\texttt{t}_{\texttt{div}}$, we employ \texttt{K-means} clustering with $k=20$, within the context of \gspinn. This clustering technique is utilized to determine the nearest cluster centroid corresponding to each pair of sample size and relative error. This approach greatly simplifies the assessment of \gspinn~and \rspinn~performance across various configurations. The resulting curves, illustrated in \Cref{Burger_qdeim_sensitivity_2}, represent \gspinn~with dashed lines and \rspinn~with solid lines. It is obvious that \gspinn~outperforms \rspinn ~almost in all the settings except for some cases that has \rspinn~has better or comparable performance which corresponds to the higher sample sizes specifically associated to the term $uu_x$.

\subsection{Korteweg-de Vries (KdV) equation}\label{subsec:KdV}
\begin{figure}[!h]
	\centering
	\includegraphics[width=1\textwidth]{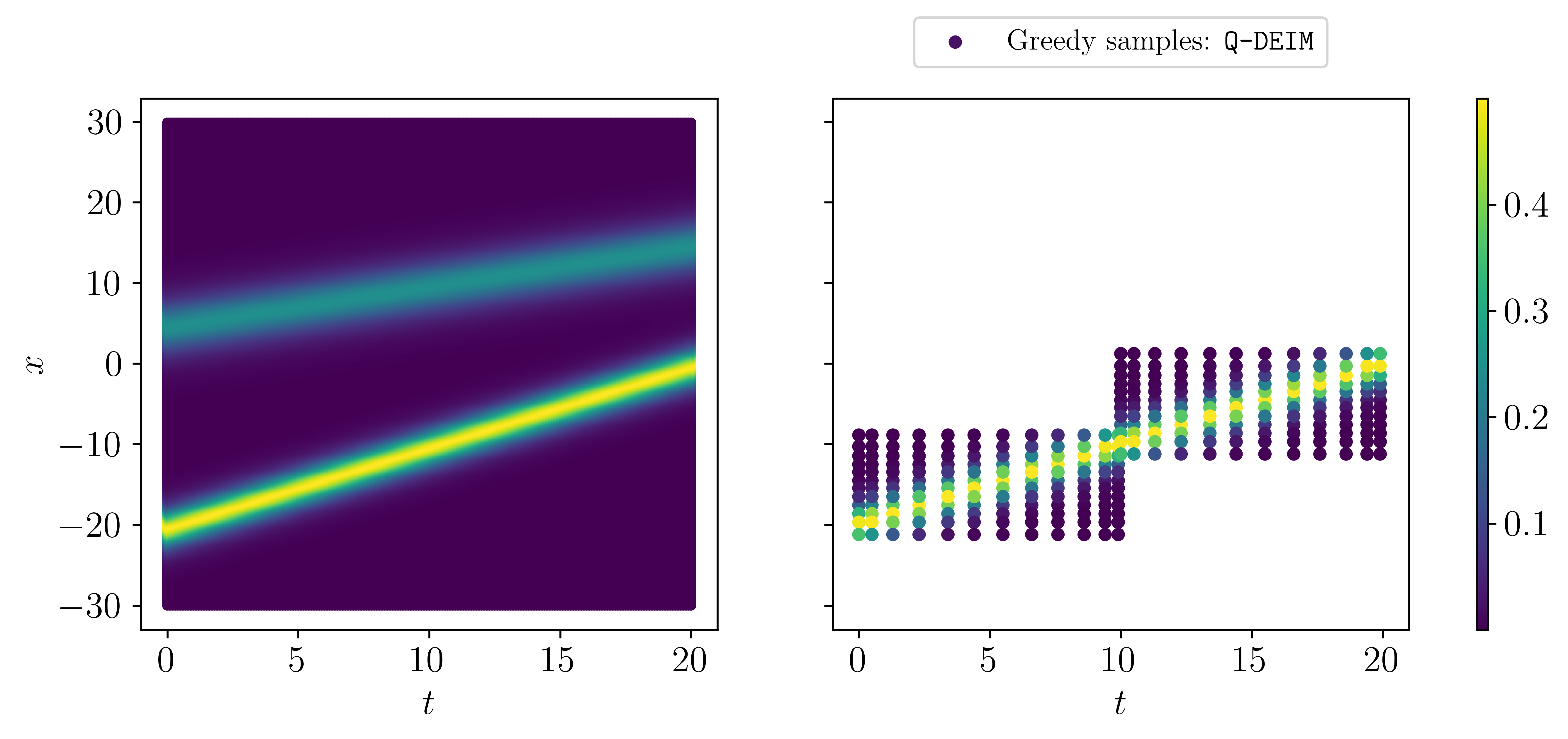}
	\caption{(left) Entire dataset ; (right) Greedy samples by \qdeim~algorithm for KdV equation  (\Cref{subsec:KdV})}
	\label{KdV_qdeim}	
\end{figure}

\begin{figure}[!h]
	\centering
	\includegraphics[width=0.9\textwidth]{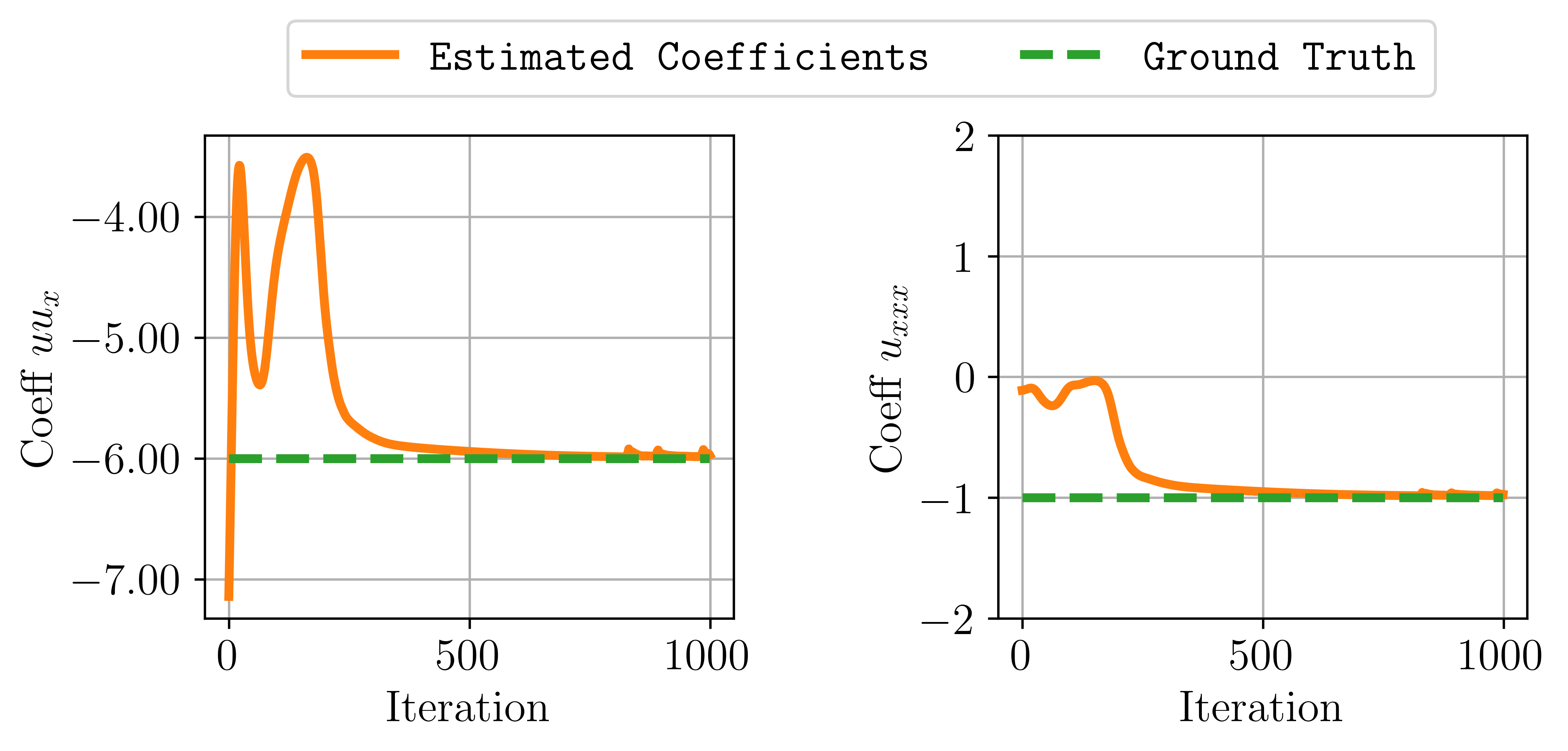}
	\caption{ Estimated coefficients by \texttt{GS-PINN} for KdV equation (\Cref{subsec:KdV})}
	\label{KdV_qdeim_coeffs}	
\end{figure}

The Korteweg-de Vries (KdV) equation is a fundamental partial differential equation in the field of mathematical physics and nonlinear waves. The KdV equation is particularly significant because it provides a mathematical description of certain types of nonlinear waves in various physical systems, including shallow water waves, plasma physics, and optics. The KdV equation can be written as a partial differential equation in the form:
\begin{equation}\label{KDV_equation}
\frac{\partial u}{\partial t } + c\ u u_x + \alpha\ u_{xxx}=0, 
\end{equation}
where $c$ and $\alpha$ are constant variables. It is known for its ability to model waves that maintain their shape while traveling through a dispersive medium, such as water with varying depths, and it can describe both solitons (solitary wave solutions) and other types of nonlinear wave phenomena. The discretization for the KdV equation, as described in \cite{rudy2017data}, is performed at $512$ space points ($n = 512$) and $201$ ($m=201$) temporal points within the original spatial domain of $x \in [-30,\ 30]$ and a time range of $t \in [0,\ 20]$, respectively. Having a snapshot matrix $\cU \in \mathbb{R}^{512 \times 201}$ demonstrates the curse of dimensionality in the training of the \pinn~for the \pde~parameter estimation, therefore we employ \qdeim~algorithm to extract the most informative samples. To do so, we set the time domain division $\texttt{t}_{\texttt{div}}=2$ and the precision value $\epsilon_{\texttt{thr}}= 10^{-3}$ whereby we earn $288$ samples which is approximately $0.28\%$ of our entire dataset. The entire dataset (left) and selected samples (right) by \qdeim~are shown in \Cref{KdV_qdeim}. Number of iteration for training is set at ${\texttt{max-iter}}=1000$ and we feed these samples into our \gspinn\ architecture to recover the parameters of the underlying \pde. The convergence of the coefficients to their true values corresponding to the terms $uu_{x}$ and $u_{xxx}$ is depicted in \Cref{KdV_qdeim_coeffs}. The computed parameter vector $p^\top = [-5.971 ,\ -0.973]$ which based on \eqref{KDV_equation} we have $c=-p_1$ and $\alpha=-p_2$. We can see that with having $0.28\%$ of the entire dataset we reach to relative errors, $\texttt{error}_1= 0.00478$ and $\texttt{error}_1= 0.0267$ corresponding to the terms $uu_x$ and $u_{xxx}$, respectively. This verifies the effectiveness of choosing most valuable and informative samples in the training of \pinn.

\begin{figure}[!h]
	\centering
	\includegraphics[width=0.8\textwidth]{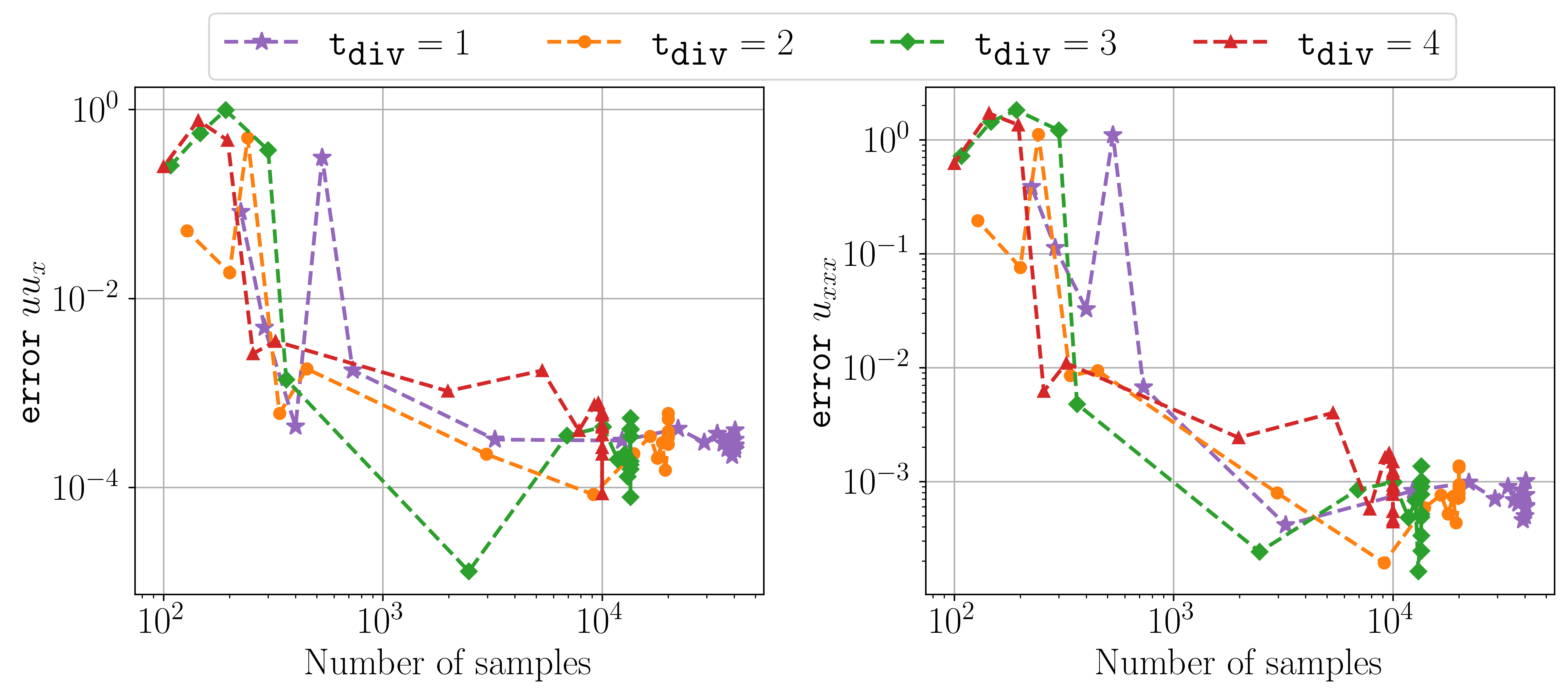}
	\caption{ Performance of \gspinn~to estimate coefficients corresponding to KdV \pde~equation (\Cref{subsec:KdV}) with different time divisions $\texttt{t}_{\texttt{div}}$ and precision values $\epsilon_{\texttt{thr}}$}
	\label{KdV_qdeim_sensitivity_1}	
\end{figure}

\Cref{KdV_qdeim_sensitivity_1} presents the relative error associated with various pairs of $( \texttt{t}_{\texttt{div}}, \epsilon_{\texttt{thr}} )$. To achieve this, we have chosen $20$ equally spaced logarithmic values for precision $\epsilon_{\texttt{thr}}$, covering the range from $10^{-10}$ to $10^{-2}$, and linked each of these values to time domain divisions $\texttt{t}_{\texttt{div}}=1,\ 2,\ 3,\ 4$. One common observation by comparing curves corresponding to the pairs $(\texttt{t}_{\texttt{div}}, \epsilon_{\texttt{thr}})$ is that for the fixed $\texttt{t}_{\texttt{div}}$ decreasing $\epsilon_{\texttt{thr}}$ has no significant impact on the number of samples that are chosen by \qdeim~and as well as corresponding relative errors. Moreover, we see that with approximately $300$ samples almost all $\texttt{t}_{\texttt{div}}$s have acceptable performance. The curve corresponding to $\texttt{t}_{\texttt{div}}=3$ with approximately $ 2500$ samples gives the minimum relative error associated to the terms $uu_x$ and $u_{xxx}$ compared to the others. Interesting observation is that time domain division $\texttt{t}_{\texttt{div}}$ limit the maximum number of samples chosen by \qdeim~at lower values of $\epsilon_{\texttt{thr}}$. 

\begin{figure}[!h]
	\centering
	\includegraphics[width=0.8\textwidth]{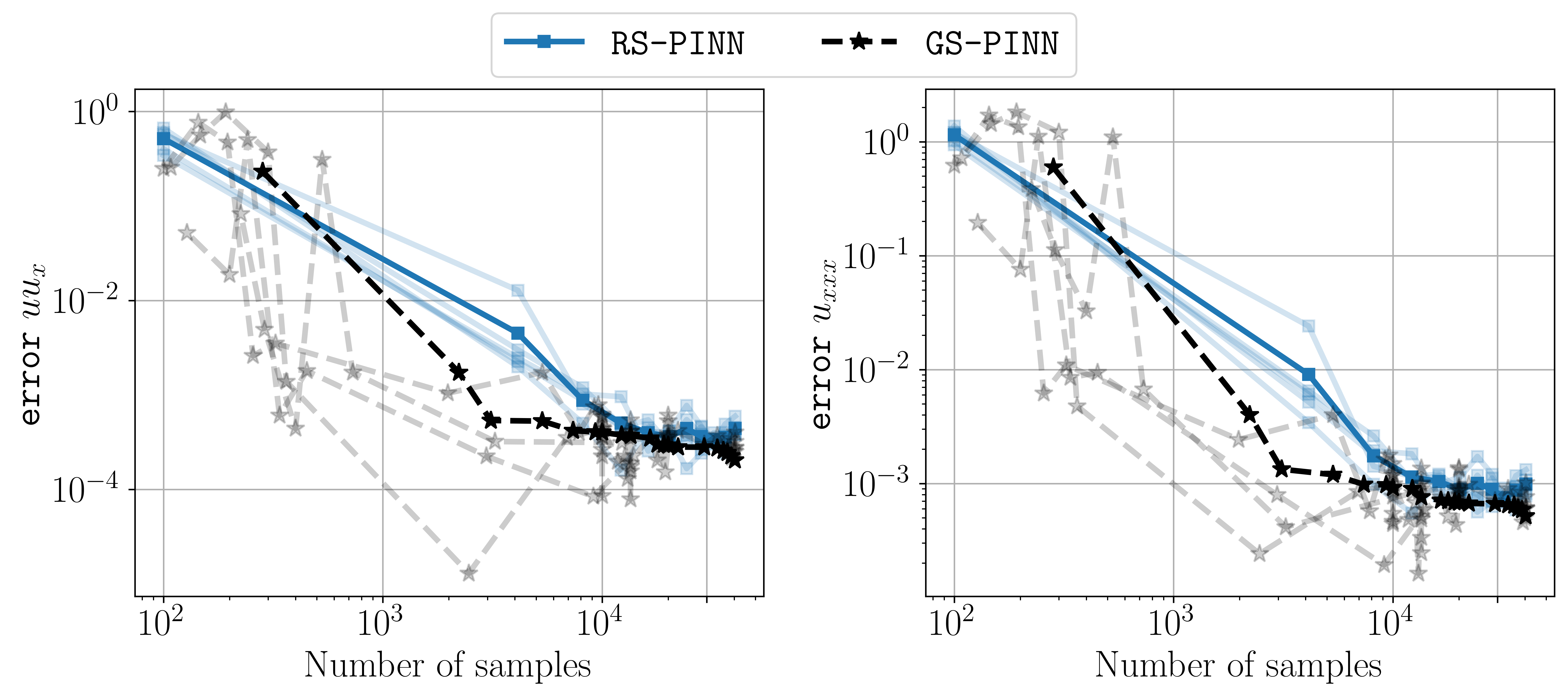}
	\caption{ Comparing \gspinn~with \rspinn~to estimate coefficients corresponding to KdV \pde~equation (\Cref{subsec:KdV}) with different time divisions $\texttt{t}_{\texttt{div}}$ and precision values $\epsilon_{\texttt{thr}}$}
	\label{KdV_qdeim_sensitivity_2}	
\end{figure}

Having the sample range through \gspinn, we are ready to evaluate the performance of \rspinn~in recovering the \pde~coefficients. The results of this comparative analysis between \gspinn~and \rspinn~are depicted in \Cref{KdV_qdeim_sensitivity_2}. To streamline the presentation of \rspinn's performance in the plots, we calculate the average relative error across five experiments for a fixed number of samples. Like the previous examples we utilize \texttt{K-means} clustering with a parameter of $k=20$ to identify the closest cluster centroid for each combination of sample size and relative error. As we can see for the lower value of sample range both \gspinn~and \rspinn~do not have good performance, however \rspinn~requires almost $4000$ samples to reach almost the same relative errors compared to \gspinn~in the sample range $200-400$ which demonstrates the significant impact of greedy sampling on \pde parameter estimation. In the higher sampling range both \rspinn~and \gspinn~have almost the same performance.


\section{Conclusion}\label{sec:conclusion}
In this paper we employed a greedy sampling approach, namely \qdeim, on the dataset resulted from descretization of the \pde s to do parameter estimation using \pinn. \qdeim~selects the most informative samples from a given snapshot matrix corresponding to a \pde. Through various simulation examples we have shown that greedy sampling approach based on \qdeim~algorithm significantly outperforms common random sampling method for the \pde~parameter estimation both in required training time for the \pinn~architecture as well as relative estimation error. Our investigation revealed that by using time domain division of the snapshot matrix of a \pde~and \qdeim~algorithm we can capture local dynamics, hence acquire more valuable samples which improves the quality of parameter estimation via \pinn. Leveraging \qdeim~and \pinn~we have estimated the corresponding parameters of three well known \pde s, (i) Allen-Cahn equation, (ii) Burgers' equation, and (iii) Korteweg-de Vries equation with $0.383\%$, $1.39\%$, and $0.28\%$ of the dataset, respectively.


\addcontentsline{toc}{section}{References}
\bibliographystyle{IEEEtran}
\bibliography{mybib}

\appendix
\section{Appendix}\label{sec:appendix}

\texttt{K-means} clustering is a popular unsupervised machine learning algorithm used for partitioning a dataset into $k$ distinct, non-overlapping subsets (clusters). The goal of the algorithm is to group similar data points together and assign them to clusters, with the number of clusters, $k$, specified by the user.

The step-by-step explanation of how the \texttt{k-means} algorithm works: (i) Initialization: Randomly select $k$ data points from the dataset as the initial centroids. The centroids are the points that will represent the center of each cluster. (ii) Assignment: assign each data point to the cluster whose centroid is closest to it. This is typically done using a distance metric, such as Euclidean distance. (iii) Update Centroids: recalculate the centroids of the clusters by taking the mean of all the data points assigned to each cluster. (iv) Repetition: repeat steps (ii) and (iii) until convergence is reached. Convergence occurs when the centroids no longer change significantly between iterations or when a certain number of iterations is reached. (v) Output: The algorithm produces $k$ clusters, and each data point is assigned to one of these clusters.

The final result of the \texttt{k-means} clustering algorithm is a set of $k$ cluster centroids and a labeling of each data point to its respective cluster. \texttt{K-means} clustering algorithm has been already implemented in Machine Learning software packages \footnote{\url{https://scikit-learn.org/stable/modules/generated/sklearn.cluster.KMeans.html}}.

To apply the \texttt{K-means} clustering algorithm on the results of the \gspinn~regarding pairs $(\texttt{t}_{\texttt{div}}, \epsilon_{\texttt{thr}})$ we construct a matrix consist of $80$ rows and $2$ columns, i.e. $ \texttt{t}_{\texttt{div}}=1,\ 2,\ 3,\ 4$ where each have corresponding $20$ relative error as mentioned in the draft. The first column is the number of samples for each pair $( \texttt{t}_{\texttt{div}}, \epsilon_{\texttt{thr}} )$ and the second column corresponds to the associated relative error. Employing \texttt{K-means} clustering algorithm with $k=20$ will result $20$ centroids each is a pair corresponds to the number of samples and relative error value. The clustering is performed with $100$ different initializations to enhance robustness.

It is worth to highlight that a \texttt{Python} package has been provided to support different implementation phase of this article which is available on \url{ https://github.com/Ali-Forootani/PINN_DEIM} or \url{https://gitlab.mpi-magdeburg.mpg.de/forootani/pinn_deim}.

\end{document}